\documentclass[11pt]{article}
\usepackage{amssymb}
\usepackage{amsfonts}
\usepackage{amsmath}
\usepackage{bm}
\usepackage{latexsym}
\usepackage{epsfig}

\newtheorem{theorem}{Theorem}[section]
\newtheorem{lemma}[theorem]{Lemma}
\newtheorem{proposition}[theorem]{Proposition}
\newtheorem{corollary}[theorem]{Corollary}

\newtheorem{definition}[theorem]{Definition}
\newtheorem{remark}[theorem]{Remark}

\newtheorem{example}[theorem]{Example}
\newenvironment{proof}{\noindent \textbf{Proof:}}{$\Box$}

\newcommand{\ignore}[1]{}

\newcommand{\enote}[1]{}
\newcommand{\knote}[1]{}
\newcommand{\rnote}[1]{}

\DeclareMathOperator{\Opt}{Opt}
\DeclareMathOperator{\Ordo}{\mathcal{O}}
\newcommand{\MaxkCSP}{\textsc{Max $k$-CSP}}
\newcommand{\MaxkCSPq}{\textsc{Max $k$-CSP$_{q}$}}
\newcommand{\MaxCSP}[1]{\textsc{Max CSP}(#1)}
\renewcommand{\Pr}{{\bf P}}
\renewcommand{\P}{{\bf P}}

\newcommand{\E}{{\bf E}}
\newcommand{\Cov}{{\bf Cov}}
\newcommand{\Var}{{\bf Var}}

\newcommand{\bits}{\{-1,1\}}

\newcommand{\Inf}{\mathrm{Inf}}

\newcommand{\eps}{\epsilon}
\newcommand{\lam}{\lambda}

\newcommand{\trunc}{\zeta}
\newcommand{\truncprod}{\chi}

\newcommand{\N}{\mathbb N}
\newcommand{\R}{\mathbb R}
\newcommand{\Z}{\mathbb Z}

\newcommand{\CalC}{{\mathcal{C}}}

\newcommand{\CalX}{{\boldsymbol{\mathcal{X}}}}
\newcommand{\CalG}{{\boldsymbol{\mathcal{G}}}}
\newcommand{\CalH}{{\boldsymbol{\mathcal{H}}}}
\newcommand{\CalY}{{\boldsymbol{\mathcal{Y}}}}
\newcommand{\CalZ}{{\boldsymbol{\mathcal{Z}}}}

\newcommand{\boldk}{{\boldsymbol k}}
\newcommand{\boldr}{{\boldsymbol r}}
\newcommand{\boldsigma}{{\boldsymbol \sigma}}

\newcommand{\hone}{{\boldsymbol{H1}}}
\newcommand{\htwo}{\boldsymbol{H2}}
\newcommand{\hthree}{\boldsymbol{H3}}
\newcommand{\hfour}{\boldsymbol{H4}}

\newcommand{\sgn}{\mathrm{sgn}}
\newcommand{\Maj}{\mathrm{Maj}}
\newcommand{\Acyc}{\mathrm{Acyc}}
\newcommand{\UniqMax}{\mathrm{UniqMax}}

\newcommand{\half}{{\textstyle \frac12}}

\newcommand{\StabThr}[2]{\Gamma_{#1}(#2)}
\newcommand{\StabThrmin}[2]{{\underline{\Gamma}}_{#1}(#2)}
\newcommand{\StabThrmax}[2]{{\overline{\Gamma}}_{#1}(#2)}
\newcommand{\TestFcn}{\Psi}

\renewcommand{\phi}{\varphi}

\begin{document}
\title{Gaussian bounds for noise correlation of functions }

\author{Elchanan Mossel\footnote{Supported by a Sloan
fellowship in Mathematics, by BSF grant 2004105, NSF Career Award (DMS 054829) and by ONR award N00014-07-1-0506. Part of this work was carried out while the author was visiting IPAM, UCLA}
\\U.C.\ Berkeley and Weizmann Institute \\mossel@stat.berkeley.edu}
\date{\today}
\maketitle

\begin{abstract}
In this paper we derive tight bounds on the expected value of
products of {\em low influence}
functions defined on correlated probability spaces.
The proofs are based on extending Fourier theory to an arbitrary number of correlated
probability spaces, on a generalization of an invariance principle
recently obtained with O'Donnell and
Oleszkiewicz for multilinear polynomials with
low influences and bounded degree and on properties of multi-dimensional
Gaussian distributions.

We present two applications of the new bounds to the theory of social choice.
We show that Majority is asymptotically the most predictable
function among all low influence functions given a random sample of the voters.
 Moreover, we derive an almost tight bound
in the context of Condorcet aggregation and low influence voting schemes on
a large number of candidates.  In particular, we show that for every low
influence aggregation function, the probability that Condorcet voting on $k$
candidates will result in a unique candidate that is preferable to all others
is $k^{-1+o(1)}$. This matches the asymptotic behavior of the majority function
for which the probability is $k^{-1-o(1)}$.

A number of applications in hardness of approximation
in theoretical computer science were obtained using the results
derived here in subsequent work by Raghavendra and by Austrin and Mossel.


A different type of applications involves hyper-graphs and arithmetic
relations in product spaces.
For example, we show that if
$A \subset \Z_m^n$ is of low influences, then the number of $k$ tuples
$(x_1,\ldots,x_k) \in A^k$ satisfying $\sum_{i=1}^k x_i \in B^n \mod m$ where
$B \subset [m]$ satisfies $|B| \geq 2$ is
$(1 \pm o(1))\P[A]^k (m^{k-1} |B|)^{n}$ which is the same as if $A$ were a random
set of probability $\P[A]$. Our results also show that for a general set $A$
without any restriction on the influences
there exists a set of coordinates $S \subset [n]$ with $|S| = O(1)$ such that
if $C = \{ x : \exists y \in A, y_{[n] \setminus S} = x_{[n] \setminus S} \}$ then
the number of $k$-tuples $(x_1,\ldots,x_k) \in C^k$ satisfying
$\sum_{i=1}^k x_i \in B^n \mod m$ is $(1 \pm o(1))\P[C]^k (m^{k-1} |B|)^{n}$.

\end{abstract}

\newpage

\section{Introduction}

\subsection{Harmonic analysis of boolean functions}
\label{sec:intro}
This paper studies low influence
functions $f : \Omega^n \to [0,1]$, where $(\Omega^n,\mu^n)$
is a product probability space and where the
\emph{influence of the $i$th coordinate on $f$},
denoted by $\Inf_i(f)$ is defined by
\begin{equation} \label{eqn:influence}
\Inf_i(f) = \E[\Var[f(X_1,\ldots,X_n) | X_j, 1 \leq j \leq n, j \neq i]],
\end{equation}
where for any set $S \subset [n]$ the conditional variance $\Var[f(X_1,\ldots,X_n) | X_i, i \in S]$ is defined via:
\[
\Var[f(X_1,\ldots,X_n) | X_i, i \in S] =
\E \Big[\big( f(X_1,\ldots,X_n) - \E[f(X_1,\ldots,X_n) | X_i, i \in S] \big)^2 |
 X_i, i \in S \Big].
 \]
The study of low influence functions is motivated by applications
from the theory of social choice in mathematical economics, by applications
in the theory of hardness of approximation in theoretical computer science and
by problems in additive number theory.
We refer the reader to some recent
papers~\cite{KKMO:04,KKMO:07,MoOdOl:05,MoOdOl:09,DiMoRe:06,SamorodnitskyTrevisan:06,Green:05} for motivation and general background.
The main theorems established here provide tight bounds on the expected
value of the product of functions defined on correlated probability spaces.
These in turn imply some new results in the theory of social choice and in
the theory of hyper-graphs. Application to hardness of approximation in
computer science were derived in subsequent work in~\cite{AustrinMossel:08} and~\cite{Raghavendra:08}.

In our main result we consider a probability measure $\P$ defined on
a space $\prod_{i=1}^k \Omega^{(i)}$.
Letting $f_i : (\Omega^{(i)})^n \to [0,1], 1 \leq i \leq k$ be a collection of
low influence functions we derive tight bounds on $\E[f_1 \ldots f_k]$
in terms of $\E[f_1],\ldots,\E[f_k]$ and a measure of correlation between the
spaces $\Omega^{(1)},\ldots,\Omega^{(k)}$.
The bounds are expressed in terms of extremal probabilities in Gaussian space, that can be calculated in the case $k=2$. When $k \geq 2$ and $\P$ is a
{\em pairwise independent distribution} our bounds show that
$\E[f_1 \ldots f_k]$ is close to $\prod_{i=1}^k \E[f_i]$.
We also apply a simple recursive argument in order to obtain results
for general functions not necessarily of low influences. The results show
that the bounds for low influence functions
hold for general functions after the
functions have been ``modified'' in a bounded number of coordinates.
The rest of the introduction is devoted to various applications
followed by statements of the main technical results.

\subsection{Prediction Of Low Influence Voting}
Suppose $n$ voters are to make a binary decision.
Assume that the outcome of the vote is determined by a {\em social choice}
function $f : \bits^n \to \bits$, so that the outcome of the vote is
 $f(x_1,\ldots,x_n)$ where $x_i \in \bits$ is the vote of voter $i$.
We assume that the votes are
independent, each $\pm 1$ with probability $\frac{1}{2}$.
It is natural to assume that the function $f$ satisfies $f(-x) = -f(x)$, i.e.,
it does not discriminate between the two candidates. Note that this implies
that $\E[f] = 0$ under the uniform distribution.
A natural way to try and
predict the outcome of the vote is to sample a subset of the voters, by sampling each voter independently with probability $\rho$.
Conditioned on a vector $X$ of votes the distribution of $Y$,
the sampled votes, is i.i.d. where $Y_i = X_i$ with probability $\rho$ and
$Y_i = \ast$ (for unknown) otherwise.

Conditioned on $Y=y$, the vector of sampled votes, the optimal prediction of
the outcome of the vote is given by $\sgn((T f)(y))$ where
\begin{equation} \label{eq:defT1}
(T f)(y) = \E[f(X) | Y = y].
\end{equation}
This implies that the probability of
correct prediction (also called predictability) is given by
\[
\P[f = \sgn(T f)] =
\frac{1}{2}(1+\E[f \, \sgn(T f)]).
\]
For example, when
$f(x) = x_1$ is the dictator function, we have
$\E[f \, \sgn(T f)] = \rho$ corresponding to the trivial fact
that the outcome of the election is known when voter $1$ is sampled
and are $\pm 1$ with probability $1/2$ otherwise.
The notion of predictability is natural in statistical contexts.
It was also studied in a more combinatorial context in~\cite{Odonnell:02}.

In the first application presented here we show that

\begin{theorem}[``Majority Is Most Predictable''] \label{thm:MIMP}
Let $0 \leq \rho \leq 1$ and $\eps > 0$ be given.  Then there
exists a $\tau > 0$ such that if $f : \bits^n \to [-1,1]$ satisfies
$\E[f] = 0$ and $\Inf_i(f) \leq \tau$ for all $i$, then
\begin{equation} \label{eq:mimp}
\E[f \, \sgn(T f) ] \leq {\textstyle \frac{2}{\pi}} \arcsin \sqrt{\rho} +
\eps,
\end{equation}
where $T$ is defined in~(\ref{eq:defT1}).
\end{theorem}

Moreover, it follows from the central limit theorem
(see Section~\ref{sec:maj};
a version of this calculation also appears in~\cite{Odonnell:02}) that
if $\Maj_n(x_1,\ldots,x_n) = \sgn(\sum_{i=1}^n x_i)$, then
\[
\lim_{n \to \infty}
\E[\Maj_n \sgn(T \Maj_n)] = {\textstyle \frac{2}{\pi}} \arcsin \sqrt{\rho}.
\]

\begin{remark}
Note that Theorem~\ref{thm:MIMP} proves a weaker statement than showing that Majority is the most predictable function. The statement only asserts that if a function has low enough influences than its predictability cannot be more than $\eps$ larger than the asymptotic predictability value achieved by the majority function when the number of voters $n \to \infty$.
This slightly inaccurate title of the theorem is inline with previous results such as the "Majority is Stablest Theorem" (see below). Similar language may be used later when informally discussing statements of various theorems.
\end{remark}

\begin{remark}
One may wonder if for a finite $n$, among {\em all} functions $f : \{-1,1\}^n \to \{-1,1\}$ with $\E[f] = 0$, majority is the most predictable function. Note that the predictability of the dictator function $f(x) = x_1$ is given by $\rho$,
and $\frac{2}{\pi}\arcsin \sqrt{\rho} > \rho$ for $\rho \to 0$. Therefore when $\rho$ is small and $n$ is large
the majority function is more predictable than the dictator function. However, note that when $\rho \to 1$
we have $\rho > \frac{2}{\pi} \arcsin \sqrt{\rho}$
and therefore for values of $\rho$ close to $1$ and large $n$ the dictator
function is more predictable than the majority function.
\end{remark}

We note that the bound obtained in Theorem~\ref{thm:MIMP} is a reminiscent of
the Majority is Stablest theorem~\cite{MoOdOl:05,MoOdOl:09} as
both involve the $\arcsin$ function. However, the two theorems are quite different. The Majority is Stablest theorem
asserts that under the same condition as in Theorem~\ref{thm:MIMP} it holds
that
\[
\E[f(X) f(Y)] \leq {\textstyle \frac{2}{\pi}} \arcsin \rho + \eps.
\]
where $(X_i,Y_i) \in \bits^2$ are i.i.d. with
$\E[X_i] = \E[Y_i] = 0$ and $\E[X_i Y_i] = \rho$. Thus ``Majority is Stablest'' considers two correlated voting vectors, while ``Majority is Most Predictable'' considers a sample of one voting vector. In fact,
both results follow from the more general invariance principle presented here.
We note a further difference between stability and predictability:
It is well known that in the context of ``Majority is Stablest'', for {\em all} $0 < \rho < 1$,
among all boolean functions with $\E[f] = 0$ the maximum of
$\E[f(x) f(y)]$ is obtained for dictator functions of the form
$f(x) = x_i$. As discussed above, for $\rho$ close to $0$ and large $n$, the dictator is less predictable than
the majority function.

We also note that the ``Ain't over until it's over'' Theorem~\cite{MoOdOl:05,MoOdOl:09} provides a bound under the same conditions on
\[
P[T f > 1-\delta],
\]
for small $\delta$.  However, this bound is not tight and does not imply Theorem~\ref{thm:MIMP}. Similarly, Theorem~\ref{thm:MIMP} does not imply the ``Ain't over until it's over'' theorem. The bounds in ``Ain't Over Until It's Over'' were
derived using invariance of $T f$ while the bound~(\ref{eq:mimp}) requires the joint invariance of $f$ and $T f$.

\subsection{Condorcet Paradoxes}
Suppose $n$ voters rank $k$ candidates. It is assumed that each voter $i$
has a linear order $\sigma_i \in S(k)$ on the candidates.
In {\em Condorcet voting},
the rankings are aggregated by deciding for each pair
of candidates which one is superior among the $n$ voters.

More formally, the aggregation results in a tournament $G_k$ on the set $[k]$.
Recall that $G_k$ is a {\em tournament} on $[k]$ if it is a
directed graph on the vertex set $[k]$ such that for all $a,b \in [k]$
either $(a > b) \in G_k$ or $(b > a) \in G_k$.
Given individual
rankings $(\sigma_i)_{i=1}^n$ the tournament $G_k$ is defined as follows.

Let $x^{a > b}(i) = 1$, if $\sigma_i(a) > \sigma_i(b)$, and $x^{a > b}(i) = -1$ if
$\sigma_i(a) < \sigma_i(b)$. Note that $x^{b > a } = -x^{a > b}$.

The binary decision between each pair of candidates is performed via
a anti-symmetric function $f : \bits^n \to \bits$ so that
$f(-x) = -f(x)$ for all $x \in \bits$.
The tournament $G_k = G_k(\sigma; f)$ is then defined
by letting $(a > b) \in G_k$ if and only if $f(x^{a > b}) = 1$.

Note that there are $2^{k \choose 2}$tournaments while there are only
$k! = 2^{\Theta(k \log k)}$ linear rankings. For the purposes of
social choice, some tournaments make more sense than others.

\begin{definition}
We say that a tournament $G_k$ is {\em linear} if it is acyclic.
We will write $\Acyc(G_k)$ for the logical statement that $G_k$ is acyclic.
Non-linear tournaments are often referred to as non-rational in economics
as they represent an order where there are $3$ candidates $a,b$ and $c$ such
that $a$ is preferred to $b$, $b$ is preferred to $c$ and
$c$ is preferred to $a$.\\

We say that the tournament $G_k$ is a {\em unique max tournament} if there is a
candidate $a \in [k]$ such that for all $b \neq a$ it holds that
$(a > b) \in G_k$.
We write $\UniqMax(G_k)$ for the logical statement that $G_k$ has a unique max.
Note that the unique max property is weaker than linearity.
It corresponds to the fact that there is a candidate that dominates all
other candidates.
\end{definition}

Following~\cite{Kalai:04,Kalai:02}, we consider the probability distribution
over $n$ voters, where the voters have independent preferences and each one chooses a ranking uniformly at random among all $k!$ orderings. Note that the marginal
distributions on vectors $x^{a > b}$ is the uniform distribution over $\bits^n$ and that if $f : \bits^n \to \bits$ is anti-symmetric then $\E[f] = 0$.\\

The case that is now understood is $k=3$.
Note that in this case $G_3$ is unique max if and only if it is linear.
Kalai~\cite{Kalai:02} studied the \emph{probability} of a rational
outcome given that the $n$ voters vote independently and at random
from the $6$ possible rational rankings.  He showed that the
probability of a rational outcome in this case may be expressed as
$\frac{3}{4}(1+\E[f Tf])$ where $T$ is the {\em Bonami-Beckner} operator with
parameter $\rho = 1/3$. The Bonami-Beckner operator may be defined as follows.
Let $(X_i,Y_i) \in \bits^2$ be i.i.d. with
$\E[X_i] = \E[Y_i] = 0$ and $\E[X_i Y_i] = \rho$ for $1 \leq i \leq n$. For $f : \bits^n \to \R$, and
$x \in \bits^n$, the Bonami-Beckner operator $T$ applied to $f$ is defined via
$(Tf)(x_1,\ldots,x_n) = \E[f(Y_1,\ldots,Y_n) | X_1 = x_1,\ldots,X_n = x_n]$.\\

It is natural to ask which function
$f$ with small influences is most likely to produce a rational
outcome.  Instead of considering small influences, Kalai
considered the essentially stronger assumption that $f$ is monotone and
``transitive-symmetric''; i.e., that for all $1 \leq i < j \leq n$
there exists a permutation $\sigma$ on $[n]$ with $\sigma(i) = j$
such that $f(x_1, \dots, x_n) = f(x_{\sigma(1)}, \dots,
x_{\sigma(n)})$ for all $(x_1, \dots, x_n)$.
Kalai conjectured that as $n \to \infty$, the maximum of
$\frac{3}{4}(1+\E[f Tf])$ among all transitive-symmetric functions approaches the same limit as
$\lim_{n \to \infty} \frac{3}{4}(1 + \E[\Maj_n \, T \Maj_n])$. This was proven using
the Majority is Stablest Theorem~~\cite{MoOdOl:05,MoOdOl:09}.
Here we obtain similar results for any value of $k$. Our result
is not tight, but almost tight. More specifically we show that:
\begin{theorem}[``Majority is best for Condorcet''] \label{thm:mibfc}
Consider Condorcet voting on $k$ candidates. Then for all $\eps > 0$ there
exists $\tau  = \tau(k,\eps) > 0$
such that if $f : \bits^n \to \bits$ is anti-symmetric
and $\Inf_i(f) \leq \tau$ for all $i$, then
\begin{equation} \label{eq:unique_max}
\P[\UniqMax(G_k(\sigma; f))] \leq k^{-1+o_k(1)} + \eps.
\end{equation}

Moreover for $f = \Maj_n$ we have $\Inf_i(f) = O(n^{-1/2})$ and it holds that
\begin{equation} \label{eq:maj_unique_max}
\P[\UniqMax(G_k(\sigma; f))] \geq k^{-1-o_k(1)} - o_n(1).
\end{equation}
\end{theorem}

Interestingly, we are not able to derive similar results for $\Acyc$.
We do calculate the probability that $\Acyc$ holds for majority.
\begin{proposition} \label{prop:maj_linear}
We have
\begin{equation} \label{eq:maj_linear}
\lim_{n \to \infty} \P[\Acyc(G_k(\sigma; \Maj_n))] =
\exp(-\Theta(k^{5/3})).
\end{equation}
\end{proposition}
We note that results in economics~\cite{Bell:81} have shown that for majority
vote the probability that the outcome will contain a Hamiltonian cycle when the
number of voters goes to infinity is $1-o_k(1)$.

\subsection{Hyper Graph and Additive Applications}
Here we discuss some applications concerning hyper-graph problems.
We let $\Omega$ be a finite set equipped with the uniform probability measure denoted $\P$.
We let $R \subset \Omega^k$ denote a $k$-wise relation.
For sets $A_1,\ldots,A_k \subset \Omega^n$ we will be interested in the
number of $k$-tuples $x^1 \in A_1,\ldots, x^k \in A_k$ satisfying the relation
$R$ in all coordinates, i.e., $(x^1_i,\ldots,x^k_i) \in R$ for all $i$.
Assume below that $R$ satisfies the following two properties:
\begin{itemize}
\item
For all $a \in \Omega$ and all $1 \leq j \leq k$ it holds that
$\P[x^i = a | (x^1,\ldots,x^k) \in R(x)] = |\Omega|^{-1}$.
(This assumption is actually not needed for the general statement -- we
state it for simplicity only).
\item
The relation $R$ is {\em connected}. This means that
for all $x,y \in R$ there exists a path $x = y(0), y(1),\ldots,y(r) = y$ in $R$ such that
$y(i)$ and $y(i+1)$ differ in one coordinate only.
\end{itemize}
We will say that the relation $R \subset \Omega^k$ is
{\em pairwise smooth} if for all
$i,j \in [k]$ and $a,b \in \Omega$ it holds that
\[
\P[x^i = a, x^j = b | (x^1,\ldots,x^k) \in R] =
\P[x^i = a | (x^1,\ldots,x^k) \in R]
\P[x^j = b | (x^1,\ldots,x^k) \in R]
\]

As a concrete example, consider the case where $\Omega = \Z_m$ and
$R$ consists of all $k$-tuples satisfying $\sum_{i=1}^k x_i \in B \mod m$ where
$B \subset \Z_m$. When $k \geq 2$ we have $\P[x_i = a | R] = m^{-1}$ for all $i$ and $a$. When $k \geq 3$, we have pairwise smoothness. The connectivity
condition holds whenever $|B| > 1$.

For a set $A \subset \Z_m^n$ and $S \subset [n]$ we define
\[
\overline{A}_S =
\{y : \exists x \in A, x_{[n] \setminus S} = y_{[n] \setminus S} \},
\underline{A}_S =
\{y : \forall x \mbox{ s.t. } x_{[n] \setminus S} = y_{[n] \setminus S}, \mbox{ it holds that } x \in A \}
\]

Our main result in the context of hyper graphs is the following.
\begin{theorem} \label{thm:ARIT}
Let $R$ be a connected relation on $\Omega^k$. Then there exist two
continuous functions $\underline{\Gamma} : (0,1)^k \to (0,1)$ and
$\overline{\Gamma} : (0,1)^k \to (0,1)$ such that
for every $\eps > 0$
there exists a $\tau > 0$ such that if $A_1,\ldots,A_k \subset \Omega^n$
are sets with
$\Inf_i(A_j) \leq \tau$ for all $i$ and $j$ then
\[
(\underline{\Gamma}(\P[A_1],\ldots,\P[A_k]) - \eps) \P[R^n]
\leq
 \P[R^n \cap (A_1,\ldots,A_k)] \leq
(\overline{\Gamma}(\P[A_1],\ldots,\P[A_k]) + \eps) \P[R^n].
\]
If $R$ is pairwise smooth, then:
\[
| \P[R^n \cap (A_1,\ldots,A_k)] - \P[R^n] \P[A_1] \ldots \P[A_k] |
\leq \eps \P[R^n].
\]
Moreover, one can take $\tau = \eps^{O(m^{k+1} \log(1/\eps)/\eps)}$.

For general sets $A_1,\ldots,A_k$, not necessarily of low influences,
there exists a set $S$ of coordinates such that
$|S| \leq O(1/\tau)$ and the statements above hold for
$\underline{A}_1^S,\ldots, \underline{A}_k^S$ and for
$\overline{A}_1^S,\ldots, \overline{A}_k^S$.
\end{theorem}

\subsection{Correlated Spaces}
A central concept that is extensively studied and repeatedly used in the
 paper is that of correlated probability spaces. The notion of correlation between two probability spaces use here
 is the same as the "maximum correlation coefficient" introduced by Hirschfeld and Gebelein~\cite{Gebelein:41}.
We will later show how to relate correlated spaces to noise operators.

\begin{definition}
Given a probability measure $\P$ defined on $\prod_{i=1}^k \Omega^{(i)}$, we
say that $\Omega^{(1)},\ldots,\Omega^{(k)}$ are {\em correlated spaces}.
For $A \subset \Omega^{(i)}$ we let
\[
\P[A] =
\P[(\omega_1,\ldots,\omega_k) \in \prod_{j=1}^k \Omega^{(j)} : \omega_i \in A],
\]
and similarly $\E[f]$ for $f : \Omega^{(i)} \to \R$. We will abuse
notation by writing $\P[A]$ for $\P^n[A]$ for
$A \subset (\prod_{i=1}^k \Omega^{(i)})^n$ or $A \subset (\Omega^{(i)})^n$ and
similarly for $\E$.

\end{definition}

\begin{definition} \label{def:rho}
Given two linear subspaces $A$ and $B$ of $L^2(\P)$ we define the
{\em correlation} between $A$ and $B$ by
\begin{equation} \label{eq:def_rho_general}
\rho(A,B; \P) = \rho(A,B) =
\sup \{\Cov[f,g] : f \in A, g \in B, \Var[f] = \Var[g] = 1\}.
\end{equation}

Let $\Omega = (\Omega^{(1)} \times \Omega^{(2)}, \P)$.
We define the {\em correlation} $\rho(\Omega^{(1)},\Omega^{(2)} ; \P)$ by letting:
\begin{equation} \label{eq:def_rho}
\rho(\Omega^{(1)},\Omega^{(2)}; \P) =
\rho(L^2(\Omega^{(1)},\P),L^2(\Omega^{(2)},\P) ; \P).
\end{equation}
More generally, let $\Omega = (\prod_{i=1}^k \Omega^{(i)}, \P)$ and
for a subset $S \subset [k]$, write $\Omega^{(S)} = \prod_{i \in S} \Omega^{(i)}$.
The {\em correlation vector}
$\underline{\rho}(\Omega^{(1)},\ldots,\Omega^{(k)} ; \P)$ is a length $k-1$
vector whose $i$'th coordinate is given by
\[
\underline{\rho}(i) = \rho(\prod_{j=1}^i \Omega^{(j)},\prod_{j=i+1}^k \Omega^{(j)} ; \P),
\]
for $1 \leq i \leq k-1$.
The {\em correlation} $\rho(\Omega^{(1)},\ldots,\Omega^{(k)} ; \P)$ is
defined
by letting:
\begin{equation} \label{eq:def_rho_k}
\rho(\Omega^{(1)},\ldots,\Omega^{(k)} ; \P) =
\max_{1 \leq i \leq k} \rho(\prod_{j = 1}^{i-1} \Omega^{(j)} \times \prod_{j=i+1}^k \Omega^{(j)},\Omega^{(i)} ; \P).
\end{equation}
\end{definition}

When the probability measure $\P$ will be clear from the context we will write
$\rho(\Omega^{(1)},\ldots,\Omega^{(k)})$ for  $\rho(\Omega^{(1)},\ldots,\Omega^{(k)} ; \P)$ etc.

\begin{remark}
It is easy to see that $\rho(\Omega^{(1)},\Omega^{(2)}; \P)$ is the
{\em second singular value} of the conditional expectation operator
mapping $f \in L^2(\Omega^{(2)},\P)$ to
$g(x) = \E[f(Y) | X = x] \in L^2(\Omega^{(1)},\P)$.
Thus $\rho(\Omega^{(1)},\Omega^{(2)}; \P)$ is the second singular value of the matrix corresponding to the operator $T$ with respect to orthonormal basis of $L^2(\Omega^{(1)},\P)$ and $L^2(\Omega^{(2)},\P)$.
\end{remark}

\begin{definition} \label{def:r_wise}
Given $(\prod_{i=1}^k \Omega^{(i)},\P)$, we say that $\Omega^{(1)},\ldots,\Omega^{(k)}$ are
{\em $r$-wise independent} if for all $S \subset [k]$ with $|S|\leq r$ and
for all $\prod_{i \in S} A_i \subset \prod_{i \in S} \Omega^{(i)}$ it holds that
\[
\P[\prod_{i \in S} A_i] = \prod_{i \in S} \P[A_i].
\]
\end{definition}
The notion of $r$-wise independence is central in computer science and discrete mathematics, in particular in the context of randomized algorithms and computational complexity.

\subsection{Gaussian Stability}
Our main result states bounds in terms of Gaussian stability measures which
we discuss next. Let $\gamma$ be the one dimensional Gaussian measure.
\begin{definition} \label{def:stabthr}
Given $\mu \in [0,1]$, define $\chi_\mu : \R \to \{0,1\}$
to be the indicator function of the interval $(-\infty, t]$, where
$t$ is chosen so that $\E_{\gamma}[\chi_\mu] = \mu$.  Explicitly, $t =
\Phi^{-1}(\mu)$, where $\Phi$ denotes the distribution function of
a standard Gaussian.  Furthermore, define
\[
\StabThrmax{\rho}{\mu,\nu} = \Pr[X \leq \Phi^{-1}(\mu)\,,\,
Y \leq \Phi^{-1}(\nu)], \quad
\]
\[
\StabThrmin{\rho}{\mu,\nu} = \Pr[X \leq \Phi^{-1}(\mu)\,,\,
Y \geq \Phi^{-1}(1-\nu)], \quad
\]
where $(X, Y)$ is a two dimensional Gaussian vector with
covariance matrix
$
\left( \begin{smallmatrix}
1 & \rho \\
\rho & 1
\end{smallmatrix} \right)
$

Given $(\rho_1,\ldots,\rho_{k-1}) \in [0,1]^{k-1}$ and
$(\mu_1,\ldots,\mu_k) \in [0,1]^k$ for $k \geq 3$
we define by induction
\[
\StabThrmax{\rho_1,\ldots,\rho_{k-1}}{\mu_1,\ldots,\mu_k} =
\StabThrmax{\rho_1}{\mu_1,\StabThrmax{\rho_2,\ldots,\rho_{k-1}}{\mu_2,\ldots,\mu_k}},
\]
and similarly $\StabThrmin{}{}$.
\end{definition}

\subsection{Statements of main results}
We now state our main results. We state the results both for low influence
functions and for general functions. For the later it is useful
to define the following notions:

\begin{definition}
Let $f : \Omega^n \to \R$ and $S \subset [n]$. We define
\[
\overline{f}^S(x) =
\sup (f(y) : y_{[n] \setminus S} = x_{[n] \setminus S}), \quad
\underline{f}^S(x) =
\inf (f(y) : y_{[n] \setminus S} = x_{[n] \setminus S}).
\]
\end{definition}

\begin{theorem} \label{thm:MIST_easy}
Let $(\prod_{j=1}^k \Omega_i^{(j)},\P_i), 1 \leq i \leq n$
be a sequence of finite probability spaces
such that for all $1 \leq i \leq n$
the minimum probability of any atom in $\prod_{j=1}^k \Omega_i^{(j)}$
is at least $\alpha$.
Assume furthermore that there exists  $\underline{\rho} \in [0,1]^{k-1}$ and $0 \leq \rho < 1$ such that
\[
\rho(\Omega^{(1)}_i,\ldots,\Omega^{(k)}_i; \P_i) \leq \rho,
\]
\begin{equation} \label{eq:rho_seq_cond}
\rho(\Omega^{(\{1,\ldots,j\})}_i,\Omega^{(\{j+1,\ldots,k\})}_i ; \P_i) \leq \underline{\rho}(j)
\end{equation}
for all $i,j$.
Then for all $\eps > 0$ there exists $\tau > 0$ such that if
\[
f_j : \prod_{i=1}^n \Omega_i^{(j)} \to [0,1].
\]
for $1 \leq j \leq k$ satisfy
\begin{equation} \label{eq:inf_bd_condition}
\max_{i,j}(\Inf_i(f_j)) \leq \tau
\end{equation}
then
\begin{equation} \label{eq:corr_bd}
\StabThrmin{\underline{\rho}}{\E[f_1],\ldots,\E[f_k]} - \eps
\leq \E[\prod_{j=1}^k f_j] \leq
\StabThrmax{\underline{\rho}}{\E[f_1],\ldots,\E[f_k]} + \eps.
\end{equation}
If we instead of (\ref{eq:rho_seq_cond}) we assume that
\begin{equation} \label{eq:rho_mat_cond}
\rho(\Omega^{(j)}_i,\Omega^{(j')}_i ; \P_i) = 0,
\end{equation}
for all $i,j \neq j'$ then
\begin{equation} \label{eq:ind_bd}
\prod_{j=1}^k \E[f_j]  - \eps
\leq \E[\prod_{j=1}^k f_j] \leq
\prod_{j=1}^k \E[f_j] + \eps.
\end{equation}
One may take
\[
\tau = \eps^{O(\frac{\log (1 / \eps) \log  (1/\alpha)}{(1-\rho) \eps})}.
\]
\end{theorem}

A truncation argument allows one to relax the conditions on the influences.

\begin{proposition} \label{prop:MIST_relaxed_easy}
For statement (\ref{eq:corr_bd}) to hold in the case where $k=2$ it suffices
to require that
\begin{equation} \label{eq:inf_bd_condition_relaxed1}
\max_{i}(\min(\Inf_i(f_1),\Inf_i(f_2))) \leq \tau
\end{equation}
instead of~(\ref{eq:inf_bd_condition}).

In the case where for each $i$ the spaces $\Omega_i^{(1)},\ldots,\Omega_i^{(k)}$
are $s$-wise independent, for statement (\ref{eq:ind_bd}) to hold it suffices
to require that for all $i$
\begin{equation} \label{eq:inf_bd_condition_relaxed2}
|\{ j : \Inf_i(f_j) > \tau\}| \leq s.
\end{equation}
\end{proposition}

An easy recursive argument allows one to conclude the following result
that does not require low influences~(\ref{eq:inf_bd_condition}).

\begin{proposition} \label{prop:MIST_easy}
Consider the setting of Theorem~\ref{thm:MIST_easy}
{\em without} the assumptions on low influences~
(\ref{eq:inf_bd_condition}).

Assuming~(\ref{eq:rho_seq_cond}), there
exists a set $S$ of size $O(1/\tau)$ such that the functions
$\overline{f}_j^S$ satisfy
\[
\E[\prod_{j=1}^k \overline{f}^S_j] \geq
\StabThrmax{\underline{\rho}}{\E[\overline{f}^S_1],\ldots,\E[\overline{f}^S_k]} - \eps \geq
\StabThrmax{\underline{\rho}}{\E[f_1],\ldots,\E[f_k]} - \eps,
\]
and the functions $\underline{f}_j^S$ satisfy
\[
\E[\prod_{j=1}^k \underline{f}^S_j] \leq
\StabThrmin{\underline{\rho}}{\E[\underline{f}^S_1],\ldots,\E[\underline{f}^S_k]} - \eps
\leq
\StabThrmin{\underline{\rho}}{\E[f_1],\ldots,\E[f_k]} - \eps.
\]
Assuming~(\ref{eq:rho_mat_cond}), we have
\[
\E[\prod_{j=1}^k \overline{f}^S_j] \geq
\prod_{j=1}^k \E[\overline{f}^S_j] - \eps \geq
\prod_{j=1}^k \E[f_j] - \eps,
\]
and similarly for $\underline{f}$.
\end{proposition}



\subsection{Road Map}
Let us review some of the main techniques we use in this paper.
\begin{itemize}

\item
We develop a Fourier theory on correlated spaces
in Section~\ref{sec:correlated}.
Previous work considered Fourier theory on one product space and
reversible operators with respect to that space~\cite{DiMoRe:06}.
Our results here allow to study non-reversible operators which in turn allows us
to study products of $k$ correlated spaces. An important fact we prove that
is used repeatedly is that general noise operators
respect ``Efron-Stein'' decomposition. This fact in particular allows us to ``truncate'' functions to their low degree parts when considering the expected
value of the product of functions on correlated spaces.

\item
In order to derive an invariance principle we need to extend the approach
of~\cite{Rotar:75,MoOdOl:05,MoOdOl:09} to prove the joint invariance of
a number of multi-linear polynomials. The proof of the extension appears
in sections~\ref{sec:background} and~\ref{sec:inv}. The proof follows the
same main steps as in~\cite{Rotar:75,MoOdOl:05,MoOdOl:09}, i.e., the Lindeberg strategy for proving the CLT~\cite{Lindeberg:22} where invariance is established by switching one variable at a time.

\item
In the Gaussian realm, we need to extend Borell's isoperimetric result~\cite{Borell:85} both in the case of two collections of Gaussians and in the case of $k > 2$ collections. This is done in Section~\ref{sec:Gaussian}.

\item
The proof of the main result, Theorem~\ref{thm:MIST_easy}
follows in Section~\ref{sec:conj}. The proof
of the extensions given in Proposition~\ref{prop:MIST_relaxed_easy}
uses a truncation
argument for which $s$-wise independence plays a crucial role. The proof
of Proposition~\ref{prop:MIST_easy} is based on a simple recursive argument.

\item
In Section~\ref{sec:social} we apply the noise bounds
in order to derive the social choice results.
Some calculations with the majority function in the social choice setting,
in particular showing the
tightness of theorems~\ref{thm:MIMP} and~\ref{thm:mibfc} are given in
Section~\ref{sec:maj}. We conclude by discussing the applications to
hyper-graphs in Section~\ref{sec:graphs}.

\end{itemize}

\subsection{Subsequent Work And Applications in Computer Science}
Subsequently to posting a draft of this paper on the Arxiv,
two applications of our results to hardness of approximation in computer
science have been established. Both results are in the context of the Unique
Games conjecture in computational complexity~\cite{Khot:02}. Furthermore,
both results consider an important problem in computer science, that is - the
problem of solving {\em constraint satisfaction problems (CSP)}.\\

Given a predicate $P : [q]^k \to \{0,1\}$, where $[q] = \{1,\ldots,q\}$ for some integer $q$,
we define $\MaxCSP{P}$ to be the algorithmic problem
where we are given a set of variables $x_1,\ldots,x_n$ taking values
in $[q]$ and a set of constraints of the form $P(l_1, \ldots, l_k)$,
where each $l_i = x_j + a$, where $x_j$ is one of the variables and $a \in [q]$
is a constant (addition is mod $q$).
More generally, in the problem of $\MaxkCSPq{}$ we are given a set of constraints
each involving $k$ of the variables $x_1,\ldots,x_n$. The most well studied case is the case of $q=2$
denoted $\MaxkCSP$.

The objective is to find an assignment to the variables satisfying as
many of the constraints as possible.
The problem of $\MaxkCSPq{}$ is NP-hard for any $k \ge 2, q \geq 2$,
and as a consequence, a large body of research is devoted to studying how well
the problem can be approximated.  We say that a (randomized)
algorithm has {\em approximation ratio} $\alpha$ if, for all
 instances, the algorithm is guaranteed to find an assignment which
 (in expectation) satisfies at least $\alpha \cdot \Opt$ of the
 constraints, where $\Opt$ is the maximum number of simultaneously
 satisfied constraints, over any assignment.\\

The results of~\cite{AustrinMossel:08} (see also~\cite{AustrinMossel:09}) 
show that assuming the Unique Games
Conjecture, for any predicate $P$ for which there exists a pairwise independent
distribution over $[q]^k$ with uniform marignals, whose support is contained
in $P^{-1}(1)$, is approximation resilient. In other words, there is no
polynomial time algorithm which achieves a better approximation factor than
assigning the variables at random. This result implies in turn that for general $k \ge 3$ and $q \ge 2$, the \MaxkCSPq{} problem is
UG-hard to approximate within $\Ordo(kq^2)/q^k + \epsilon$. Moreover,
for the special case of $q = 2$, i.e., boolean variables, it gives hardness
of $(k + \Ordo(k^{0.525}))/2^k + \epsilon$, improving
upon the best previous bound~\cite{SamorodnitskyTrevisan:06} of $2k/{2^k} + \epsilon$ by essentially a factor $2$.
Finally, again for $q=2$, assuming that the famous Hadamard
Conjecture is true, the results are further improved, and the
$\Ordo(k^{0.525})$ term can be replaced by the constant $4$.

These results should be compared to prior work
by Samordnitsky and Trevisan~\cite{SamorodnitskyTrevisan:06}
who using the Gowers norm, proved that the
\MaxkCSP{} problem has a hardness factor of $2^{\lceil \log_2 k+1
 \rceil}/2^k$, which is $(k+1)/2^k$ for $k = 2^r-1$, but can be as
 large as $2k/2^k$ for general $k$.

 From the quantitative point of view~\cite{AustrinMossel:09} gives stronger
stronger hardness than~\cite{SamorodnitskyTrevisan:06} for
  $\MaxkCSPq{}$, even in the already thoroughly explored $q = 2$ case.
  These improvements may seem very small, being an improvement only by
  a multiplicative factor $2$.  However, it is well known that it is
  impossible to get non-approximability results which are better than
  $(k+1)/2^k$, and thus, in this
  respect, the hardness of $(k+4)/2^k$ assuming the Hadamard
  Conjecture is in fact optimal to within a very small \emph{additive}
  factor.  Also, the results of~\cite{AustrinMossel:09}
  give approximation resistance of
  $\MaxCSP{P}$ for a much larger variety of predicates (any $P$
  containing a balanced pairwise independent distribution).

  From a qualitative point of view, the analysis of~\cite{AustrinMossel:09}
  is very direct. Furthermore, it is
  general enough to accommodate any domain $[q]$ with virtually no
  extra effort.  Also, their proof using the main result of the current paper,
  i.e.,  bounds on expectations of
  products under certain types of correlation, putting it in the same
  general framework as many other UGC-based hardness results, in
  particular those for $2$-CSPs.\\

  In a second beautiful result by Raghavendra~\cite{Raghavendra:08}
  the results of the current paper were used to obtain very general
  hardness results for $\MaxCSP{P}$. In~\cite{Raghavendra:08} it is
  shown that for every predicate $P$ and for every approximation factor
  which is smaller than the UG-hardness of the problem, there exists a
  polynomial time algorithm which achieves this approximation ratio.
  Thus for every $P$ the UG-hardness of $\MaxCSP{P}$ is sharp.
  The proof of the results uses the results obtained here in order to
  define and analyze the reduction from UG given the integrality gap
  of the corresponding convex optimization problem. We note that for most
  predicates the UG hardness of $\MaxCSP{P}$ is unknown and therefore the
  results of~\cite{Raghavendra:08}
  complement those of~\cite{AustrinMossel:08}.

\subsection{Acknowledgments}
I would like to thank
Noam Nisan for suggesting that generalization of the invariance
principle should be useful in the social choice context and
Gil Kalai for pointing out some helpful references.
I would like to thank Terence Tao for helpful
discussions and references on additive
number theory and Szemeredi regularity.
Thanks to Per Austrin, Nick Crawford, Marcus Issacson and
Ryan O'Donnell for a careful reading of a draft of this paper.
Finally, many thanks to an anonymous referee for many helpful comments.

\section{Correlated Spaces and Noise} \label{sec:correlated}

In this section we define and study the notion of correlated spaces and
noise operators in a general setting.

\subsection{Correlated Probability Spaces and Noise Operators}
We begin by defining noise operators and giving some basic examples.
\begin{definition} \label{def:markov_operator}
Let $(\Omega^{(1)} \times \Omega^{(2)},\P)$ be two correlated spaces.
The {\em Markov Operator} associated with $(\Omega^{(1)},\Omega^{(2)})$ is the
operator mapping $f \in L^p(\Omega^{(2)},\P)$ to $Tf \in L^p(\Omega^{(1)},\P)$ by:
\[
(T f)(x) = \E[f(Y) | X = x],
\]
for $x \in \Omega^{(1)}$ and where $(X,Y) \in \Omega^{(1)} \times \Omega^{(2)}$ is distributed according to $\P$.
\end{definition}

\begin{example}
In order to define {\em Bonami-Beckner} operator $T = T_{\rho}$ on a space
$(\Omega,\mu)$, consider the space $(\Omega \times \Omega,\nu)$ where
$\nu(x,y) = (1-\rho) \mu(x) \mu(y) + \rho \delta(x=y) \mu(x)$, where $\delta(x=y)$ is the function on $\Omega \times \Omega$ which takes the value $1$ when $x=y$, and $0$ otherwise.
In this case, the operator $T$ satisfies:
\begin{equation} \label{eq:bonami_beckner1}
(Tf)(x) = \E[f(Y) | X = x],
\end{equation}
where the conditional distribution of $Y$ given $X = x$ is
$\rho \delta_x + (1-\rho) \mu$, where $\delta_x$ is the delta measure on $x$.
\end{example}

\begin{remark}
The construction above may be generalized as follows. Given any Markov chain on
$\Omega$ that is reversible with respect to $\mu$, we may look at the measure
$\nu$ on $\Omega \times \Omega$ defined by the Markov chain. In this case $T$
is the Markov operator determined by the chain. The same construction applies
under the weaker condition that $T$ has $\mu$ as its stationary distribution.
\end{remark}

It is straightforward to verify that:
\begin{proposition} \label{prop:T_general}
Suppose that for each $1 \leq i \leq n$, $(\Omega_i^{(1)} \times \Omega_i^{(2)},\mu_i)$
are correlated spaces and
$T_i$ is the Markov operator associated with $\Omega_i^{(1)}$ and $\Omega_i^{(2)}$.
Then $(\prod_{i=1}^n \Omega_i^{(1)},\prod_{i=1}^n \Omega_i^{(2)},\prod_{i=1}^n \mu_i)$
defines two correlated spaces  and the Markov operator $T$ associated with
them is given by $T = \otimes_{i=1}^n T_i$.
\end{proposition}

\begin{example} \label{ex:bonami_beckner}
For product spaces $(\prod_{i=1}^n \Omega_i, \prod_{i=1}^n \mu_i)$, the
Bonami-Beckner operator $T = T_{\rho}$ is defined by
\begin{equation} \label{eq:bonami_beckner2}
T = \otimes_{i=1}^n T^i_{\rho},
\end{equation}
where $T^i$ is the Bonami-Beckner operator on
$(\Omega_i \times \Omega_i, \mu_i)$.
This Markov operator is the one most commonly discussed in previous work,
see e.g.~~\cite{KaKaLi:88,KKMO:04,MoOdOl:09}.
In a more recent work~\cite{DiMoRe:06} the case of $\Omega_i \times \Omega_i$
with $T_i$ a reversible Markov operator with respect to a measure $\mu_i$
on $\Omega_i$ was studied.
\end{example}

\begin{example} \label{ex:sample}
In the context of the Majority is most Predictable Theorem~\ref{thm:MIMP}, the underlying space is
$\Omega = \{\pm 1\} \times \{0, \pm1\}$ where element $(x,y) \in \Omega$
corresponds to a voter with vote $x$ and a sampled vote $y$ where either
$y=x$ if the vote is queried or $y=0$ otherwise.
The probability measure $\mu$ is given by
\[
\mu(x,y) = \frac{1}{2} (\delta(x=1) + \delta(x=-1))
(\rho \delta(y = x) + (1-\rho) \delta(y=0)).
\]

Note that the marginal distributions on $\Omega_S = \{0,\pm 1\}$
and $\Omega_V = \{\pm 1\}$ are given by
\[
\mu = (1-\rho) \delta_{0} + \frac{\rho}{2} (\delta_{-1} + \delta_1),\,\,\,
\nu = \frac{1}{2} (\delta_{-1} + \delta_1),
\]
and
\[
\nu(\cdot | \pm 1) = \delta_{\pm 1},\,\,\,
\nu(\cdot | 0)  = \frac{1}{2}(\delta_1 + \delta_{-1}).
\]
Given independent copies $\mu_i$ of $\mu$ and $\nu_i$ of $\nu$,
the measure $\mu = \otimes_{i=1}^n \mu_i$
corresponds to the distribution of a sample of voters where
each voter is sampled independently with probability $\rho$ and the
distribution of the voters is given by $\nu = \otimes_{i=1}^n \nu_i$.
\end{example}

\begin{example} \label{ex:condorcet}
The second non-reversible example is natural in the context
of Condorcet voting. For simplicity, we first discuss the case of
$3$ possible outcomes. The general case is discussed later.

Let $\tau$ denote the uniform measure on the set permutations on the set $[3]$ denoted $S_3$. Note that each element
$\sigma \in S_3$ defines an element $f \in \{-1,1\}^{3 \choose 2}$ by letting
$f(i,j) = \sgn(\sigma(i) - \sigma(j))$. The measure so defined, defines $3$ correlated probability spaces $(\{\pm 1\}^{3 \choose 2},\P)$.

Note that the projection of $\P$ to each coordinate is uniform and
\[
\P(f(3,1) = -1 | f(1,2) = f(2,3) = 1) = 0,\,\,\,
\P(f(3,1) = 1 | f(1,2) = f(2,3) = -1) = 0,\,\,\,
\]
and
\[
\P(f(3,1) = \pm 1 | f(1,2) \neq f(2,3)) = 1/2.
\]
\end{example}

\subsection{Properties of Correlated Spaces and Markov Operators}
Here we derive properties of correlated spaces and Markov operators that will
be repeatedly used below. We start with the following which was already known to R{\'e}nyi~\cite{Renyi:59}.
\begin{lemma} \label{lem:rho_cond_exp}
Let $(\Omega^{(1)} \times \Omega^{(2)},\P)$ be two correlated spaces.
Let $f$ be a $\Omega^{(2)}$ measurable function with $\E[f] = 0$, and
$\E[f^2] = 1$. Then among all $g$ that are $\Omega^{(1)}$
measurable satisfying $\E[g^2] = 1$, a maximizer of $|\E[fg]|$ is given by
\begin{equation} \label{eq:rho_cond_exp}
g = \frac{Tf}{\sqrt{\E[(Tf)^2]}},
\end{equation}
where $T$ is the Markov operator associated with $(\Omega^{(1)},\Omega^{(2)})$.
Moreover,
\begin{equation} \label{eq:rho_cond_exp2}
|\E[g f]| = \frac{|\E[f Tf]|}{\sqrt{\E[(Tf)^2]}} = \sqrt{\E[(Tf)^2]}.
\end{equation}
\end{lemma}

\begin{proof}
To prove (\ref{eq:rho_cond_exp})
let $h$ be an $\Omega^{(1)}$ measurable function with $\| h \|_2 = 1$.
Write $h = \alpha g + \beta h'$ where $\alpha^2 + \beta^2 = 1$ and
$\| h' \|_2 = 1$ is orthogonal to $g$.
From the properties of conditional
expectation it follows that $\E[f h'] = 0$. Therefore we may choose an
optimizer satisfying $\alpha \in \pm 1$. Equation~(\ref{eq:rho_cond_exp2}) follows since $Tf$ is a conditional expectation. The same reasoning shows that
$\E[fg] = 0$ for every $\Omega^{(1)}$ measurable function $g$ if $Tf$ is identically $0$.\\
\end{proof}

The following lemma is useful in bounding $\rho(\Omega^{(1)},\Omega^{(2)} ; \P)$ from Definition~\ref{def:rho}
in generic situations. Roughly speaking, it shows that connectivity of the support of
$\P$ on correlated spaces $\Omega^{(1)} \times \Omega^{(2)}$ implies that $\rho < 1$.
\begin{lemma} \label{lem:cheeger}
Let $(\Omega^{(1)} \times \Omega^{(2)},\P)$ be two correlated spaces
such that the probability of the smallest atom in $\Omega^{(1)} \times \Omega^{(2)}$ is
at least $\alpha > 0$. Define a bi-partite graph
$G = (\Omega^{(1)},\Omega^{(2)},E)$
where $(a,b) \in \Omega^{(1)} \times \Omega^{(2)}$ satisfies $(a,b) \in E$ if
$\P(a,b) > 0$.
Then if $G$ is connected then
\[
\rho(\Omega^{(1)},\Omega^{(2)} ; \P) \leq 1 - \alpha^2/2.
\]
\end{lemma}

\begin{proof}
For the proof it would be useful to consider
$G' = (\Omega^{(1)} \cup \Omega^{(2)},E)$, a weighted
{\em directed} graph where the weight $W(a,b)$ of the directed edge from $a$ to $b$
is $\P[b | a]$ and the weight of the directed edge from $b$ to $a$ is $W(b,a) = \P[a | b]$.
 Note that the minimal
non-zero weight must be at least $\alpha$ and that
$W(a \to b) > 0$ iff $W(b > a) > 0$. This later fact implies that $G'$
is strongly connected. Note furthermore that $G'$ is bi-partite.

Let $A$ be the transition probability matrix defined by the weighted graph
$G'$. Since $G'$ is connected and $W(a,b) \geq \alpha$ for all $a$ and
$b$ such that $W(a \to b) > 0$,
it follows by Cheeger's inequality
that the spectral gap of $A$ is at least $\alpha^2/2$.
Since $G'$ is connected and bi-partite, the multiplicities of the eigenvalues $\pm 1$ are both $1$.
Corresponding eigenfunctions are the constant $1$ functions and the function taking the value $1$ on $\Omega^{(1)}$
and the value $-1$ on $\Omega^{(2)}$.

In order to bound $\rho$ it suffices by Lemma~\ref{lem:rho_cond_exp} to bound $\| Af \|_2$ for a function $f$ that is supported on $\Omega^{(2)}$ and satisfies $\E[f] = 0$.  Note that such a function is orthogonal to the eigen-vectors of $A$ corresponding to the eigenvalues $-1$ and $1$. It therefore
follows that $\| A f \|_2 \leq (1-\alpha^2/2) \| f \|_2$ as needed.
\end{proof}

One nice property of Markov operators that will 
be used below is that they respect the Efron-Stein decomposition.
Given a vector $x$ in an $n$ dimensional product space and $S \subset [n]$ we
write $x_S$ for the vector $(x_i : i \in S)$. Given probability spaces
$\Omega_1,\ldots,\Omega_n$, we use the convention of writing $X_i$ for a random variable that is distributed according to the measure of $\Omega_i$ and $x_i$ for an element of $\Omega_i$. We will also write $X_S$ for $(X_i : i \in S)$.

\begin{definition} \label{def:efron_stein}
Let $(\Omega_1,\mu_1),\ldots,(\Omega_n,\mu_n)$ be discrete probability spaces
$(\Omega,\mu) = \prod_{i=1}^n (\Omega_i,\mu_i)$. The Efron-Stein decomposition of
$f : \Omega \to \R$ is given by
\begin{equation} \label{eq:efronstein}
f(x) = \sum_{S \subseteq [n]} f_S(x_S),
\end{equation}
where the functions $f_S$ satisfy:
\begin{itemize}
\item
$f_S$ depends only on $x_S$.
\item
For all $S \not \subseteq S'$ and all $x_{S'}$ it holds that:
\[
\E[f_S | X_{S'} = x_{S'}] = 0.
\]
\end{itemize}

\end{definition}
It is well known that the Efron-Stein decomposition exists and
that it is unique~\cite{EfronStein:81}. We quickly recall the proof of existence.
The function $f_S$ is given by:
\[
f_S(x) = \sum_{S' \subseteq S} (-1)^{|S \setminus S'|} \E[f(X) | X_{S'} = x_{S'}]
\]
which implies
\[
\sum_{S} f_S(x) = \sum_{S'} \E[f | X_{S'} = x] \sum_{S : S' \subseteq S} (-1)^{|S \setminus S'|} =
\E[f | X_{[n]} = x_{[n]}] = f(x).
\]
Moreover, for $S \not \subseteq S'$ we have $\E[f_S | X_{S'} = x_{S'}] = \E[f_S | X_{S' \cap S} = x_{S' \cap S}]$
and for $S'$ that is a strict subset of $S$ we have:
\begin{eqnarray*}
\E[f_S | X_{S'} = x_{S'}] &=& \sum_{S'' \subset S} (-1)^{|S \setminus S''|} \E[f(X) |X_{S'' \cap S'} = x_{S'' \cap S'}] \\ &=& \sum_{S'' \subset S'} \E[f(X) |X_{S''} = x_{S''}] \sum_{S''  \subset \tilde{S} \subset S'' \cup (S \setminus S')} (-1)^{|S \setminus \tilde{S}|} = 0.
\end{eqnarray*}

We now prove that the Efron-Stein decomposition ``commutes'' with Markov
operators.
\begin{proposition} \label{prop:TEfronStein1}
Let $(\Omega_i^{(1)} \times \Omega_i^{(2)},\P_i)$ be correlated spaces and let
$T_i$ the Markov operator associated with $\Omega_i^{(1)}$ and $\Omega_i^{(2)}$ for
$1 \leq i \leq n$. Let
\[
\Omega^{(1)} = \prod_{i=1}^n \Omega_i^{(1)},\quad
\Omega^{(2)} = \prod_{i=1}^n \Omega_i^{(2)},\quad
\P = \prod_{i=1}^n \P_i, \quad
T = \otimes_{i=1}^n T_i.
\]
 Suppose $f \in L^2(\Omega^{(2)})$ has Efron-Stein decomposition~(\ref{eq:efronstein}). Then the Efron-Stein decomposition of $Tf$ satisfies:
\[
(Tf)_S = T (f_S).
\]
\end{proposition}

\begin{proof}
Clearly $T (f_S)$ is a function of $x_S$ only.
Moreover, for all $S \not \subseteq S'$ and $x^0_{S'}$ it holds that
\[
\E[(T (f_S))(X) | X_{S'} = x^0_{S'}] =
\E[f_S(Y) | X_{S'} = x^{0}_{S'}] =
\E[ \E[f_S(Y) | Y_{S'}] \P[Y_{S'} | X_{S'} = x^{0}_{S'}]] = 0,
\]
where the second equality follows from the fact that $Y$ is independent of $X_{S'}$ condition on
$Y_{S'}$.
\end{proof}

We next derive a useful bound showing that in the setting above if
$\rho(\Omega_i^{(1)} \times \Omega_i^{(2)}; \P) < 1$ for all $i$ then $Tf$ depends on
the ``low degree expansion'' of $f$.
\begin{proposition} \label{prop:TEfronStein2}
Assume the setting of Proposition~\ref{prop:TEfronStein1} and that further
for all $i$ it holds that $\rho(\Omega_i^{(1)},\Omega_i^{(2)} ; \P_i) \leq \rho_i$.
Then for all $f$ it holds that
\[
\| T (f_S) \|_2 \leq \left( \prod_{i \in S} \rho_i \right) \| f_S \|_2.
\]
\end{proposition}

\begin{proof}
Without loss of generality if suffices to prove the statement of Proposition~\ref{prop:TEfronStein2}
for $S = [n]$. Thus our goal is to snow that
\[
\| T f \|_2 \leq \left( \prod_{i=1}^n \rho_i \right) \| f \|_2.
\]
(from now on $S$ will denote a set different than $[n]$).
For each $0 \leq r \leq n$, let $T^{(r)}$ denote the following operator.
$T^{(r)}$ maps a function $g$ of
$z = (z_1,\ldots,z_n) = (x_1,\ldots,x_{r-1},y_r,\ldots,y_n)$
to a function $T^{(r)} g$ of
$w = (w_1,\ldots,w_n) = (x_1,\ldots,x_r,y_{r+1},\ldots,y_n)$
defined as follows:
\[
T^{(r)} g(w) =
\E[g(Z) | W = w].
\]
(Here $Z = (X_1,\ldots,X_{r-1},Y_r,\ldots,Y_n)$
and similarly  $W$).

Let $g$ be a function such that for any subset $S \subsetneq [n]$
and all $z_S$,
\[
\E[g(Z) | Z_S = z_S] = 0.
\]
We claim that
\begin{equation} \label{eq:cond_exp_contr}
\| T^{(r)}g \|_2 \leq \rho_r \| g \|_2
\end{equation}
and that for all subsets
$S \subsetneq [n]$ it holds that
\begin{equation} \label{eq:cond_exp_0}
 \E[(T^{(r)}g)(W) | W_S = w_S] = 0.
\end{equation}
Note that~(\ref{eq:cond_exp_contr}) and~(\ref{eq:cond_exp_0})
together imply the desired bound as
$T = T^{(n)} \cdots T^{(1)}$.

For~(\ref{eq:cond_exp_contr}) note that if $S = [n] \setminus \{r\}$  and
$f = T^{(r)}g$ then by lemma~\ref{lem:rho_cond_exp}
\[
\E[f^2(W) | Z_S = z_S] =
\E[g(Z) f(W) | Z_S = z_S] \leq
\rho_r \sqrt{\E[f^2(W) | Z_S = z_S] \E[g^2(Z) | Z_S = z_S]}.
\]
So
\[
\E[f^2(W) | Z_S = z_S] \leq \rho_r^2 \E[g^2(W) | Z_S = z_S],
\]
which gives $\| f \|_2 \leq \rho_r \| g \|_2$ by integration.

For~(\ref{eq:cond_exp_0}) we note that if $S \subsetneq [n]$ then
\[
\E[f(W) | W_S = w_S] = \E[g(Z) | W_S = w_S] =
\E[ \E[g(Z) | Z_S] \P[Z_S| W_S = w_S]] = 0.
\]
This concludes the proof.
\end{proof}

\begin{proposition} \label{prop:TEfronStein3}
Assume the setting of Proposition~\ref{prop:TEfronStein1}.
Then
\[
\rho(\prod_{i=1}^n \Omega_i^{(1)},\prod_{i=1}^n \Omega_i^{(2)} ; \prod_{i=1}^n \P_i) =
\max \rho(\Omega_i^{(1)}, \Omega_i^{(2)}).
\]
\end{proposition}

\begin{proof}
Let $f \in L^2(\prod_{i=1}^n \Omega_i^{(2)})$ with $\E[f] = 0$ and $\Var[f] = 1$.
Expand $f$ according to its Efron-Stein decomposition
\[
f = \sum_{S \neq \emptyset} f_S,
\]
($f_{\emptyset} = 0$ since $\E[f] = 0$).
Then by propositions~\ref{prop:TEfronStein1} and~\ref{prop:TEfronStein2}
\begin{eqnarray*}
\E[ (Tf)^2] &=& \E[(T (\sum_S f_S))^2] = \E[(\sum_S T f_S)^2] =
\sum_S \E[(T f_s)^2] \leq
\sum_{S \neq \emptyset} \prod_{i \in S} \rho_i \| f_S \|_2^2 \\ &\leq&
\max_i \rho_i^2 \sum_{S \neq \emptyset} \| f_S \|_2^2 = \max_i \rho_i^2.
\end{eqnarray*}
The other inequality is trivial.
\end{proof}



\section{Background: Influences and Hypercontractivity} \label{sec:background}
In this section we recall and generalize some
definitions and results from~\cite{MoOdOl:09}.
In particular, the generalizations allow us the study
non-reversible Markov operators and correlated ensembles.
For the reader who is familiar with~\cite{MoOdOl:09} it suffices to look
at subsections~\ref{subsec:vec_multi_lin} and~\ref{subsec:vec_hyper}.

\subsection{Influences and noise stability in product spaces} \label{sec:general}

Let $(\Omega_1, \mu_1), \dots, (\Omega_n, \mu_n)$ be
probability spaces and let $(\Omega, \mu)$ denote the
product probability space.  Let
\[
f = (f_1,\ldots,f_k): \Omega_1 \times \cdots \times \Omega_n \to \R^k.
\]

\begin{definition} \label{def:influence_general}
The \emph{influence of the $i$th coordinate on $f$} is
\[
\Inf_i(f) = \sum_{1 \leq j \leq k} \E [\Var[f_j | x_1,\ldots,x_{i-1},x_{i+1},\ldots,x_n]],
\]
\end{definition}

\subsection{Multi-linear Polynomials} \label{sec:setup}

In this sub-section we recall and slightly generalize
the setup and notation used in~\cite{MoOdOl:09}.
Recall that we are interested in
functions on product of finite probability spaces, $f : \Omega_1
\times \cdots \times \Omega_n \to \R$.  For each $i$, the space of
all functions $\Omega_i \to \R$ can be expressed as the span of a
finite set of orthonormal random variables, $X_{i,0} = 1,
X_{i,1}, X_{i, 2}, X_{i,3}, \dots$; then $f$ can be
written as a multilinear polynomial in the $X_{i,j}$'s. In
fact, it will be convenient for us to mostly disregard the
$\Omega_i$'s and work directly with sets of orthonormal random
variables; in this case, we can even drop the restriction of
finiteness.  We thus begin with the following definition:

\begin{definition}\cite{MoOdOl:09}
We call a collection of finitely many orthonormal real random
variables, one of which is the constant $1$, an \emph{orthonormal
ensemble}.  We will write a typical \emph{sequence} of $n$
orthonormal ensembles as $\CalX = (\CalX_1, \dots, \CalX_n)$,
where $\CalX_i = \{X_{i,0} = 1, X_{i,1}, \dots,
X_{i,m_i}\}$.  We call a sequence of orthonormal ensembles
$\CalX$ \emph{independent} if the ensembles are independent families
of random variables.

We will henceforth be concerned only with independent sequences of
orthonormal ensembles, and we will call these \emph{sequences of
ensembles}, for brevity. Similarly, when writing an ensemble we will always mean an
orthogonal ensemble.
\end{definition}

\begin{remark}  \label{rem:simple}\cite{MoOdOl:09} Given a sequence of independent random variables
$X_1, \dots, X_n$ with $\E[X_i] = 0$ and
$\E[X_i^2] = 1$, we can view
them as a sequence of ensembles $\CalX$ by renaming $X_i =
X_{i,1}$ and setting $X_{i,0} = 1$ as required.
\end{remark}

\begin{definition}\cite{MoOdOl:09}  We denote by $\CalG$ the \emph{Gaussian sequence of
ensembles}, in which $\CalG_i = \{G_{i,0} = 1, G_{i,1},
G_{i,2}, \dots\,G_{i,m_i}\}$ and all $G_{i,j}$'s with $j \geq 1$
are independent standard Gaussians.
\end{definition}
The Gaussian ensembles discussed in this paper will often have $m_i$ chosen to match the $m_i$ of a given
ensemble. \\

As mentioned, we will be interested in \emph{multilinear
polynomials} over sequences of ensembles.  By this we mean sums of
products of the random variables, where each product is obtained
by multiplying one random variable from each ensemble.
\begin{definition}\cite{MoOdOl:09}
A \emph{multi-index} $\boldsigma$ is a sequence $(\sigma_1, \dots,
\sigma_n)$ in $\N^n$. The \emph{degree} of $\boldsigma$, denoted
$|\boldsigma|$, is $|\{i \in [n] : \sigma_i > 0\}|$.  Given a
doubly-indexed set of indeterminates $\{x_{i,j}\}_{i \in [n], j
\in \N}$, we write $x_\boldsigma$ for the monomial $\prod_{i =
1}^n x_{i,\sigma_i}$.  We now define a \emph{multilinear
polynomial} over such a set of indeterminates to be any expression
 \begin{equation} \label{eqn:Q}
Q(x) = \sum_{\boldsigma} c_\boldsigma x_\boldsigma
 \end{equation}
where the $c_\boldsigma$'s are real constants, all but finitely
many of which are zero. The \emph{degree} of $Q(x)$ is
$\max\{|\boldsigma| : c_\boldsigma \neq 0\}$, at most $n$.  We
also use the notation
\[
Q^{\leq d}(x) = \sum_{|\boldsigma| \leq d} c_\boldsigma
x_\boldsigma
\]
and, analogously, $Q^{= d}(x)$ and $Q^{> d}(x)$.
\end{definition}

Naturally, we will consider applying multilinear polynomials $Q$
to sequences of ensembles $\CalX$; the distribution of these
random variables $Q(\CalX)$ is the subject of our invariance
principle.  Since $Q(\CalX)$ can be thought of as a function on a
product space $\Omega_1 \times \cdots \times \Omega_n$ as
described at the beginning of this section, there is a consistent
way to define the notions of influences, $T_\rho$, and noise
stability from Section~\ref{sec:general}.  For example, the
``influence of the $i$th ensemble on $Q$'' is
\[
\Inf_i(Q(\CalX)) = \E[\Var[Q(\CalX) \mid \CalX_1, \dots,
\CalX_{i-1}, \CalX_{i+1}, \dots, \CalX_n]].
\]
Using independence and orthonormality, it is easy to show the
following formulas:
\begin{proposition} \label{prop:infQ}
Let $\CalX$ be a sequence of ensembles and $Q$ a multilinear
polynomial as in~(\ref{eqn:Q}). Then
\[
 \E[Q(\CalX)] = c_{\bf{0}}; \qquad \E[Q(\CalX)^2] = \sum_{\boldsigma}
 c_\boldsigma^2; \qquad
 \Var[Q(\CalX)] = \sum_{|\boldsigma| > 0} c_\boldsigma^2;
\]
and
\[
 \Inf_i(Q(\CalX)) = \sum_{\boldsigma : \sigma_i > 0} c_\boldsigma^2;
\]
\end{proposition}



For $\rho \in [0,1]$ we define the operator $T_\rho$ as acting
formally on multilinear polynomials $Q(x)$ as in~(\ref{eqn:Q}) by
\begin{equation} \label{eq:bonami_beckner3}
(T_\rho Q)(x) = \sum_{\boldsigma} \rho^{|\boldsigma|} c_\boldsigma
x_\boldsigma.
\end{equation}
We note the definition in~(\ref{eq:bonami_beckner1}) and~(\ref{eq:bonami_beckner2}) are consistent with the definition in~(\ref{eq:bonami_beckner3}) in the sense that for any ensemble $\CalX$ the two definitions results in the
same function $(T_{\rho} Q)(\CalX)$.

\bigskip
We finally recall the notion of ``low-degree
influences'', a notion that has proven crucial in the analysis of
PCPs in hardness of approximation in computer science
(see, e.g., \cite{KKMO:04}).
\begin{definition}\cite{MoOdOl:09}
\label{def:low-degree-influence} The \emph{$d$-low-degree influence of the $i$th ensemble on $Q(\CalX)$}
is
\[
\Inf^{\leq d}_i(Q(\CalX)) = \Inf^{\leq d}_i(Q) = \sum_{\boldsigma:
|\boldsigma| \leq d, \sigma_i > 0} c_\boldsigma^2.
\]
Note that this gives a way to define low-degree influences
$\Inf^{\leq d}_i(f)$ for functions $f : \Omega_1 \times \cdots
\Omega_n \to \R$ on finite product spaces.
\end{definition}
There isn't an especially natural interpretation of $\Inf_i^{\leq
d}(f)$.  However, the notion is important for PCPs due to the fact
that a function with variance $1$ cannot have too many coordinates
with substantial low-degree influence; this is reflected in the
following easy proposition:
\begin{proposition} \label{prop:infD}\cite{MoOdOl:09}
Suppose $Q$ is multilinear polynomial as in~(\ref{eqn:Q}).  Then
\[
\sum_i \Inf_i^{\leq d}(Q) \leq d \cdot \Var[Q].
\]
\end{proposition}
The proof follows since:
\[
\sum_i \Inf_i^{\leq d}(Q) = \sum_i \sum_{\boldsigma:
|\boldsigma| \leq d, \sigma_i > 0} c_\boldsigma^2 =
\sum_{\boldsigma: 0 < |\boldsigma| \leq d} |\boldsigma| c_\boldsigma^2 \leq
d \sum_{\boldsigma: 0 < |\boldsigma|} c_\boldsigma^2 = d \Var[Q].
\]

\subsection{Vector valued multi-linear polynomials} \label{subsec:vec_multi_lin}
For the invariance principle discussed here
we will need to consider vector-valued multi-linear polynomials.
\begin{definition}
A $k$-dimensional \emph{multilinear
polynomial} over a set of indeterminates is given by
\begin{equation} \label{eqn:QD}
Q = (Q_1,\ldots,Q_k)
\end{equation}
where each $Q_j$ is a multi-linear polynomial as in~(\ref{eqn:Q}).
The \emph{degree} of $Q$ is the maximal degree of the $Q_j$'s.
\end{definition}

\begin{definition}
We adopt the standard notation and write $\| Q \|_q$ for
$\E^{1/q}[\sum_{j=1}^k |Q_i|^q]$; we write $\Var[Q]$ for $\E[\| Q - \E[Q] \|_2^2]$
and
\[
\Inf_i(Q(\CalX)) = \sum_{j=1}^k \Inf_i(Q_j(\CalX))
\]
\end{definition}

Using these definitions, it is easy to see that
\begin{proposition} \label{prop:infQDIM}
Let $\CalX$ be a sequence of ensembles and $Q = (Q_1,\ldots,Q_k)$ a $k$
dimensional multilinear polynomial where $Q_j$ is defined as in~(\ref{eqn:Q}) with
$c^j_{\boldsigma}$ as its coefficients. Then:
\[
 \E[Q(\CalX)] = (c^1_{\bf{0}},\ldots,c^k_{\bf{0}});
\qquad \| Q(\CalX) \|_2^2 = \sum_{j,\boldsigma} {c^j}_\boldsigma^2; \qquad
 \Var[Q(\CalX)] = \sum_{j,|\boldsigma^j| > 0} {c^j}_\boldsigma^2.
\]
\end{proposition}

Finally, we recall the standard multi-index notation associated
with $k$-dimensional multi-linear polynomials. A \emph{multi-index} $\bf{i}$
of dimension $k$ is a vector $(i_1,\ldots,i_k)$, where each
$i_j$ is an integer. We write $|\bf{i}|$ for
$i_1 + \cdots + i_k$ and $\bf{i}!$ for $i_1 ! i_2 ! \cdots i_k!$.
Given a function $\psi$ of $k$ variables, we will write $\psi^{(\bf{i})}$ for the partial derivative of $f$ taken $i_1$ times with respect to the first variable, $i_2$ with respect to the second etc. (we will only consider functions $\psi$ that are smooth enough that the order of derivatives does not matter). We will also write $Q^{\bf{i}}$ for the product
$Q_1^{i_1}  \cdots  Q_k^{i_k}$.

\subsection{Hypercontractivity}
As in~\cite{MoOdOl:09} the invariance
principle requires that the ensembles involved are hypercontractive.
Recall that $Y$ is $(2,q,\eta)$-hypercontractive with some $\eta \in (0,1)$ if and
only if $\E[Y]=0$ and $\E[|Y|^{q}]<\infty.$ Also, if
$Y$ is $(2,q,\eta)$-hypercontractive then $\eta \leq
(q-1)^{-1/2}.$\\


\begin{definition}
Let $\CalX$ be a sequence of ensembles.  For $1 \leq p \leq q <
\infty$ and $0 < \eta < 1$ we say that $\CalX$ is
\emph{$(p,q,\eta)$-hypercontractive} if
\[
\|(T_\eta Q)(\CalX)\|_q \leq \|Q(\CalX)\|_p
\]
for every multilinear polynomial $Q$ over $\CalX$.
\end{definition}
Since $T_{\eta}$ is a contractive semi-group, we have
\begin{remark} If $\CalX$ is $(p,q,\eta)$-hypercontractive then it
is $(p,q,\eta')$-hypercontractive for any $0 < \eta' \leq \eta$.
\end{remark}

There is a related notion of hypercontractivity for \emph{sets} of
random variables which considers all polynomials in the variables,
not just multilinear polynomials; see, e.g.,
Janson~\cite{Janson:97}. We summarize some of the basic properties below,
see~\cite{MoOdOl:09} for details.

\begin{proposition}~\cite{MoOdOl:09} \label{prop:join-hypercon}
Suppose $\CalX$ is a sequence of $n_1$ ensembles and $\CalY$ is an
independent sequence of $n_2$ ensembles.  Assume both are
$(p,q,\eta)$-hypercontractive.  Then the sequence of ensembles
$\CalX \cup \CalY = (\CalX_1, \dots, \CalX_{n_1}, \CalY_1, \dots,
\CalY_{n_2})$ is also $(p,q,\eta)$-hypercontractive.
\end{proposition}

\begin{proposition}~\cite{MoOdOl:09} \label{prop:hypercon}
Let $\CalX$ be a $(2,q,\eta)$-hypercontractive sequence of
ensembles and $Q$ a multilinear polynomial over $\CalX$ of degree
$d$.  Then
\[
\|Q(\CalX)\|_q \leq \eta^{-d} \; \|Q(\CalX)\|_2.
\]
\end{proposition}

We end this section by recording the some hypercontractive
estimates to be used later. The result for $\pm 1$
Rademacher variables is well known and due originally to
Bonami~\cite{Bonami:70} and independently
Beckner~\cite{Beckner:75}; the same result for Gaussian and uniform
random variables is also well known and in fact follows easily from the
Rademacher case. The optimal hypercontractivity constants for general
finite spaces was recently determined by Wolff~\cite{Wolff:07}
(see also~\cite{Oleszkiewicz:03}):
\begin{theorem} \label{thm:bonami} Let $X$ denote either a uniformly random $\pm 1$
bit, a standard one-dimensional Gaussian, or a random variable
uniform on $[-\sqrt{3}, \sqrt{3}]$. Then $X$ is $(2, q,
(q-1)^{-1/2})$-hypercontractive.
\end{theorem}
\begin{theorem} \label{thm:wolff}~\cite{Wolff:07} (Wolff)\ \ Let $X$ be any
mean-zero random variable on a finite probability space in which
the minimum nonzero probability of any atom is $\alpha \leq 1/2$.
Then $X$ is $(2, q, \eta_q(\alpha))$-hypercontractive, where
\[
\eta_q(\alpha) = \left(\frac{A^{1/q'} - A^{-1/q'}}{A^{1/q} -
A^{-1/q}} \right)^{-1/2}
\]
\[
\text{with } \quad A = \frac{1-\alpha}{\alpha}, \quad 1/q + 1/q' =
1.
\]
\end{theorem}
Note the following special case:
\begin{proposition} \label{prop:wolff}~\cite{Wolff:07}
\[
\eta_3(\alpha) = \left(A^{1/3} + A^{-1/3}\right)^{-1/2} \;\;
\mathop{\sim}^{\alpha \to 0} \quad \alpha^{1/6},
\]
and also
\[
\half \alpha^{1/6} \leq \eta_3(\alpha) \leq 2^{-1/2},
\]
for all $\alpha \in [0,1/2]$.
\end{proposition}


\subsection{Vector Hyper-Contraction} \label{subsec:vec_hyper}
For our purposes we will also need to obtain hypercontraction results in cases where $Q$ is a $k$-dimensional multi-linear polynomial.
We will need to consider vector-valued multi-linear polynomials.
\begin{proposition} \label{prop:hyperconmulti}
Let $\CalX$ be a $(2,q,\eta)$-hypercontractive sequence of
ensembles and $Q$ a multilinear polynomial over $\CalX$ of degree
$d$ and dimension $k$. Assume $q$ is integer and let $\bf{i}$ be a multi-index
with $|{\bf{i}}| = q$. Then
\[
\E[|Q^{\bf{i}}(\CalX)|] \leq \eta^{-dq} \prod_{j=1}^k \E[Q_j^2]^{i_j/2},
\]
\end{proposition}

\begin{proof}
\[
\E[|Q^{\bf{i}}(\CalX)|] \leq \prod_{j=1}^k \E[|Q_j|^q]^{i_j/q}
\leq \eta^{-dq} \prod_{j=1}^k \E[Q_j^2]^{i_j/2},
\]
where the first inequality is H\"{o}lder and the second follows by
hypercontractivity.
\end{proof}

\section{Multi-dimensional Invariance principle} \label{sec:inv}
In this section we generalize the invariance principle from~\cite{MoOdOl:09} to
the multi-dimensional setting. We omit some easy steps that are either
identical or easy adaptation of the proofs of~\cite{MoOdOl:09}.

\subsection{Hypotheses for invariance theorems}
Below we will prove a generalization of the
invariance principle~\cite{MoOdOl:09}. The invariance principle proven there
concerns a multilinear polynomial $Q$ over two
hypercontractive sequences of ensembles, $\CalX$ and $\CalY$;
furthermore, $\CalX$ and $\CalY$ are assumed to satisfy a
``matching moments'' condition, described below.

It is possible to generalize the invariance principle to vector valued
multi-linear polynomials under each of the hypercontractivity assumptions
$\hone,\htwo,\hthree$ and $\hfour$ of~\cite{MoOdOl:09}.
However, since the proof of all generalizations is essentially the same and
since for the applications studied here it suffices to consider the
hypothesis $\hthree$, this is the only hypothesis that will be
discussed in the paper. It is defined as follows:

\begin{enumerate}

\item[$\hthree$] Let $\CalX$ be a sequence of $n$
ensembles in which the random variables in each ensemble $\CalX_i$
form a basis for the real-valued functions on some finite
probability space $\Omega_i$.  Further assume that the least
nonzero probability of any atom in any $\Omega_i$ is $\alpha \leq
1/2$, and let $\eta = \half \alpha^{1/6}$.  Let $\CalY$ be any
$(2,3,\eta)$-hypercontractive sequence of ensembles such that all the random
variables in $\CalY$ are independent.
Finally, let $Q$ be a $k$ dimensional multilinear polynomial as
in~(\ref{eqn:QD}).\\

\end{enumerate}

\subsection{Functional Setting}
The essence of our invariance principle is that if $Q$ is of
bounded degree and has low influences then the random variables
$Q(\CalX)$ and $Q(\CalY)$ are close in distribution. The simplest
way to formulate this conclusion is to say that if $\TestFcn : \R^k
\to \R$ is a sufficiently nice ``test function'' then
$\TestFcn(Q(\CalX))$ and $\TestFcn(Q(\CalY))$ are close in
expectation.

\begin{theorem} \label{thm:cj}  Assume hypothesis $\hthree$.
Further assume that $Q$ is a $k$ dimensional multi-linear polynomial,
that $\Var[Q] \leq 1$, $\deg(Q) \leq d$,
and $\Inf_i(Q) \leq \tau$ for all $i$.
Let $\TestFcn : \R^k \to \R$ be a $\CalC^3$
function with $|\TestFcn^{(\bf{i})}| \leq B$ uniformly for every
vector $\bf{i}$ with $|\bf{i}|=3$.  Then
\[
\Bigl| \E\bigl[\TestFcn(Q(\CalX))\bigr] -
\E\bigl[\TestFcn(Q(\CalY))\bigr] \Bigr| \leq
\eps := 2 d B\,k^3 (8\alpha^{-1/2})^d\;\tau^{1/2}.
\]
\end{theorem}

\begin{proof}
Note that by Proposition~\ref{prop:wolff}, the random variables satisfy $(2,3,\eta)$
hypercontractivity with $\eta = \frac{1}{2} \alpha^{1/6}$.

We begin by defining intermediate sequences between $\CalX$ and
$\CalY$.  For $i = 0, 1, \dots, n$, let $\CalZ^{(i)}$ denote the
sequence of $n$ ensembles $(\CalY_1, \dots, \CalY_i, \CalX_{i+1},
\dots, \CalX_n)$ and let $Q^{(i)} = Q(\CalZ^{(i)})$.  Our
goal will be to show
\begin{equation} \label{eqn:bound}
\Bigl|\E\bigl[\TestFcn(Q^{(i-1)})\bigr] -
\E\bigl[\TestFcn(Q^{(i)})\bigr]\Bigr| \leq
2B k^3\,\eta^{-3d} \Inf_i(Q)^{3/2}
\end{equation}
for each $i \in [n]$.  Summing this over $i$ will complete the
proof since $\CalZ^{(0)} = \CalX$, $\CalZ^{(n)} = \CalY$, and
\[
\sum_{i=1}^n \Inf_i(Q)^{3/2} \leq \tau^{1/2} \cdot \sum_{i=1}^n \Inf_i(Q) =
\tau^{1/2} \cdot \sum_{i=1}^n \Inf_i^{\leq d}(Q)
\leq d \tau^{1/2},
\]
where we used Proposition~\ref{prop:infD} and $\sum_j \Var[Q_j] \leq 1$.

\newcommand{\Qt}{\tilde{Q}} \newcommand{\Rt}{R} \newcommand{\St}{S}
Let us fix a particular $i \in [n]$ and proceed to
prove~(\ref{eqn:bound}).  Given a multi-index $\boldsigma$, write
$\boldsigma \setminus i$ for the same multi-index except with
$\sigma_i = 0$.  Now write
\begin{eqnarray*}
\Qt &=& \sum_{\boldsigma : \sigma_i = 0} c_\boldsigma
\CalZ^{(i)}_\boldsigma, \\
\Rt &=& \sum_{\boldsigma : \sigma_i > 0} c_\boldsigma
X_{i,\sigma_i} \cdot \CalZ^{(i)}_{\boldsigma \setminus i}, \\
\St &=& \sum_{\boldsigma : \sigma_i > 0} c_\boldsigma
Y_{i,\sigma_i} \cdot \CalZ^{(i)}_{\boldsigma \setminus i}.
\end{eqnarray*}
Note that $\Qt$ and the variables $\CalZ^{(i)}_{\boldsigma
\setminus
  i}$ are independent of the variables in $\CalX_i$ and $\CalY_i$ and
that $Q^{(i-1)} = \Qt + \Rt$ and $Q^{(i)} = \Qt + \St$.
\\

To bound the left side of~(\ref{eqn:bound}) --- i.e.,
$|\E[\TestFcn(\Qt + \Rt) - \TestFcn(\Qt + \St)]|$ --- we use
Taylor's theorem: for all $x, y \in \R$,
\[
\Bigl|\TestFcn(x+y) - \sum_{|\boldk| < 3}
\frac{\TestFcn^{\boldk}(x)\;y^{\boldk}}{\boldk!}\Bigr| \leq \sum_{|\boldk| = 3}
\frac{B}{{\boldk}!}\,|y|^{\boldk}.
\]
In particular,
\begin{equation} \label{eq:R_taylor}
\Bigl|\E[\TestFcn(\Qt + \Rt)] - \sum_{|\boldk| < 3}
\E\Bigl[\frac{\TestFcn^{(\boldk)}(\Qt)\;\Rt^{\boldk}}{\boldk!}\Bigr]\Bigr| \leq
\sum_{|\boldk| = 3} \frac{B}{\boldk!}\,\E\bigl[|\Rt|^{\boldk}\bigr]
\end{equation}
and similarly,
\begin{equation} \label{eq:S_taylor}
\Bigl|\E[\TestFcn(\Qt + \St)] - \sum_{|\boldk| < 3}
\E\Bigl[\frac{\TestFcn^{(\boldk)}(\Qt)\;\St^{\boldk}}{\boldk!}\Bigr]\Bigr| \leq
\sum_{|\boldk| = 3} \frac{B}{\boldk!}\,\E\bigl[|\St|^{\boldk}\bigr]
\end{equation}
We will see below that $\Rt$ and $\St$ have finite $3$'rd moments.
Moreover,
for $0 \leq \boldk \leq \boldr$ with $|\boldr| = 3$ it holds that
$|\TestFcn^{(\boldk)}(\Qt)\,\Rt^{\boldk}| \leq
|\boldk!\,B\,\Qt^{\boldr-\boldk}\,\Rt^{\boldk}|$ (and similarly for $\St$). Thus all
moments above are finite.
We now claim
that for all $0 \leq |\boldk| < 3$ it holds that
\begin{equation} \label{eq:S_T_moments}
\E[\TestFcn^{(\boldk)}(\Qt)\,\Rt^{\boldk}] = \E[\TestFcn^{(\boldk)}(\Qt)\,\St^{\boldk}].
\end{equation}
This follows from the fact that the expressions in the expected values
when viewed as multi-linear polynomials in the variables in $\CalX_i$ and
$\CalY_i$ respectively are of degree $\leq 2$ and each monomial term in $\CalX_i$ has the same coefficient as the corresponding monomial in $\CalY_i$.


From (\ref{eq:R_taylor}), (\ref{eq:S_taylor}) and
(\ref{eq:S_T_moments}) it follows that
\begin{equation} \label{eq:taylor_bd}
 |\E[\TestFcn(\Qt + \Rt) - \TestFcn(\Qt + \St)]| \leq \sum_{\boldr = 3}
\frac{B}{\boldr!}\,(\E[|\Rt|^{\boldr}] + \E[|\St|^{\boldr}]).
\end{equation}
We now use hypercontractivity.  By
Proposition~\ref{prop:join-hypercon} each $\CalZ^{(i)}$ is
$(2,r,\eta)$-hypercontractive.  Thus by
Proposition~\ref{prop:hyperconmulti},
\begin{equation} \label{eq:R_S_hyper}
\E[|\Rt|^{\boldr}] \leq \eta^{-3d} \prod_{j=1}^k \E[\Rt_j^2]^{r_j/2}, \quad
\E[|\St|^{\boldr}] \leq \eta^{-3d} \prod_{j=1}^k \E[\St_j^2]^{r_j/2}.
\end{equation}
However,
\begin{equation} \label{eq:S_R_inf}
\E[\St_j^2] = \E[\Rt_j^2] = \sum_{\boldsigma^j : \sigma^j_i > 0}
c_{\boldsigma^j}^2 = \Inf_i(Q_j) \leq \Inf_i(Q).
\end{equation}
Combining (\ref{eq:taylor_bd}), (\ref{eq:R_S_hyper}) and
(\ref{eq:S_R_inf}) it follows that
\[
 |\E[\TestFcn(\Qt + \Rt) - \TestFcn(\Qt + \St)]| \leq
2 B k^3\,\eta^{-3d} \cdot\Inf_i(Q)^{3/2}
\]
confirming~(\ref{eqn:bound}) and completing the proof.
\end{proof}

\subsection{Invariance principle --- other functionals, and smoothed version}

The basic invariance principle shows that $\E[\TestFcn(Q(\CalX))]$
and $\E[\TestFcn(Q(\CalY))]$ are close if $\TestFcn$ is a
$\CalC^3$ functional with bounded $3$rd derivative.  To show that
the distributions of $Q(\CalX)$ and $Q(\CalY)$ are close in other
senses we need the invariance principle for less smooth
functionals.  This we can obtain using straightforward
approximation arguments, see for example~\cite{MoOdOl:09}.
For applications involving bounded functions, it will be important to bound
the following functionals.
We let $f_{[0,1]} : \R \to \R$ be defined by
\[
f_{[0,1]}(x) = \max(\min(x,1),0) = \min(\max(x,0),1),
\]
and
$\trunc : \R^k \to \R$ be defined by
\begin{equation} \label{eq:def_trunc}
\trunc(x) = \sum_{i=1}^k (x_i - f_{[0,1]}(x_i))^2.
\end{equation}
Similarly, we define
\begin{equation} \label{eq:def_truncprod}
\truncprod(x) = \prod_{i=1}^k f_{[0,1]}(x_i).
\end{equation}

Repeating the proofs of~\cite{MoOdOl:09} one obtains:
\begin{theorem} \label{thm:smooththeorem}
Assume hypothesis $\hthree$. Further assume that $Q = (Q_1,\ldots,Q_k)$
is a $k$-dimensional multi-linear polynomial with
$\Var[Q] \leq 1$ and $\Inf_i^{\leq\,\log(1/\tau)/K}(Q) \leq \tau$
for all $i$, where
\[
 K = \log(1/\alpha).
\]
Suppose further that for all $d$ it holds that
$\Var[Q^{> d}] \leq (1-\gamma)^{2d}$ where $0 < \gamma < 1$.
Write $R = Q(\CalX)$
and $S = Q(\CalY)$. Then
\[
 \Bigl|\E\bigl[\trunc(R)\bigr] -
 \E\bigl[\trunc(S)\bigr]\Bigr| \leq
 \tau^{\Omega(\gamma/K)},
\]
where $\Omega(.)$ hides a constant depending only on $k$.
Similarly,
\[
 \Bigl|\E\bigl[\truncprod(R)\bigr] -
 \E\bigl[\truncprod(S)\bigr]\Bigr| \leq
 \tau^{\Omega(\gamma/K)}.
\]
\end{theorem}

\begin{proof}
The proof for $\trunc$ uses the fact that
the function $\trunc$ admits
approximations $\trunc_{\lam}$ such that
\begin{itemize}
\item
\[
\| \trunc - \trunc_{\lam} \|_{\infty} = O(\lam^2),
\]
\item
\[
\|\trunc_{\lam}^{\bf{r}} \|_{\infty} = O(\lam^{-1}),
\]
for all $\bf{r}$ with $|\bf{r}| = 3$.
\end{itemize}
This implies that for all $k$ dimensional degree $d$ polynomials we have:
\begin{equation} \label{eq:apx1}
 \Bigl|\E\bigl[\trunc(R)\bigr] -
 \E\bigl[\trunc(S)\bigr]\Bigr| \leq O(\alpha^{-d/3} \tau^{-1/3}).
\end{equation}
See~\cite{MoOdOl:09} for details. In order to obtain the result for polynomials
with decaying tails we use the fact that
\[
|\trunc(a_1 + b_1,\ldots,a_k + b_k) - \trunc(a_1,\ldots,a_k)|
\leq \sum_{i=1}^k (|a_i b_i| + b_i^2).
\]
This implies that evaluating $\trunc$ at polynomials truncated at level $d$ results in an
error of at most $O(\exp(-d \gamma))$ which together with the bound~(\ref{eq:apx1})implies the desired bound for $\trunc$.

The proof for $\truncprod$ is similar as the function $\trunc$ admits
approximations $\truncprod_{\lam}$ such that
\begin{itemize}
\item
\[
\| \truncprod - \truncprod_{\lam} \|_{\infty} = O(\lam),
\]
\item
\[
\|\truncprod_{\lam}^{\bf{r}} \|_{\infty} = O(\lam^{-2}),
\]
for all $\bf{r}$ with $|\bf{r}| = 3$.
\end{itemize}
\end{proof}

Of particular interest to us is the following corollary.

\begin{corollary} \label{cor:smooththeorem}
Assume hypothesis $\hthree$.
For each $1 \leq i \leq n$ and
$1 \leq j \leq k$, let $\CalX_i^j$ denote an ensemble all of whose
elements are linear combinations of elements in $\CalX_i$.
Similarly for each $1 \leq i \leq n$ and
$1 \leq j \leq k$, let $\CalY_i^j$ denote an ensemble all of whose
elements are linear combinations of elements in $\CalY_i$.

Assume further that for all $i$ and $j$
it holds that $|\CalY_i^j| = |{\CalX}_i^j|$
and that there is a one to one correspondence between $Y \in \CalY_i^j$
given by:
\[
Y = \sum_{r=0}^{|\CalY_i|-1} a(i,j,k) Y_k,
\]
and $X \in \CalX_i^j$ given by:
\[
X = \sum_{r=0}^{|\CalX_i|-1} a(i,j,k) X_k,
\]
where $\CalY_i = \{Y_0,\ldots,Y_{|\CalY_i|-1}\}$ and
$\CalX_i = \{X_0,\ldots,X_{|\CalY_i|-1}\}$.

Let $Q = (Q_1,\ldots,Q_k)$ a be multi-linear polynomial with
$\Var[Q] \leq 1$ and $\Inf_i^{\leq\,\log(1/\tau)/K}(Q) \leq \tau$ for all $i$
where
\[
 K = \log(1/\alpha).
\]
Suppose further that for all $d$ it holds that
$\Var[Q_j^{> d}] \leq (1-\gamma)^{2d}$ where $0 < \gamma < 1$.
Write $R = (Q_1(\CalX^1),\ldots,Q_k(\CalX^k))$
and $S = (Q_1(\CalY^1),\ldots,Q_k(\CalY^k))$. Then
\[
 \Bigl| \E\bigl[\trunc(R)\bigr]- \E\bigl[\trunc(S)\bigr]\Bigr| \leq
 \tau^{\Omega(\gamma/K)},
\]
and
\[
 \Bigl| \E\bigl[\truncprod(R)\bigr]-
 \E\bigl[\truncprod(S)\bigr]\Bigr| \leq
 \tau^{\Omega(\gamma/K)},
\]
where the $\Omega(\cdot)$ hides a constant depending only on $k$.
\end{corollary}

\begin{proof}
The proof follows immediately from the previous theorem noting that
$\Inf_i$ and $\Var[Q^{> d}]$ are basis independent.
\end{proof}

\section{Noise in Gaussian Space} \label{sec:Gaussian}
In this section we derive the Gaussian bounds correlation bounds needed
for our applications. The first bound derived in subsection~\ref{subsec:borell}  is an easy extension of~\cite{Borell:85}.
The second one gives a quantitative estimate on iterations of the first bound
that will be needed for some of the applications.
\subsection{Noise stability in Gaussian space} \label{subsec:borell}
We begin by recalling some definitions and results relevant for
``Gaussian noise stability''.  Throughout this section we consider
$\R^n$ to have the standard $n$-dimensional Gaussian distribution,
and our probabilities and expectations are over this
distribution. Recall Definition~\ref{def:stabthr}.
We denote by
$U_\rho$ the Ornstein-Uhlenbeck operator acting on
$L^2(\R^n)$ by
\[
(U_\rho f)(x) = \E[f(\rho x + \sqrt{1-\rho^2}\, y)],
\]
where $y$ is a standard $n$-dimensional Gaussian.

It is immediate to see that if $f,g \in L^2(\R^n)$ then
\[
\E[f U_{\rho} g] = \E[g U_{\rho} f] =
\E[f(X_1,\ldots,X_n) g(Y_1,\ldots,Y_n)]
\]
where $(X_i, Y_i)$ are independent two dimensional Gaussian vectors with
covariance matrix
$
\left( \begin{smallmatrix}
1 & \rho \\
\rho & 1
\end{smallmatrix} \right)
$
The results of Borell~\cite{Borell:85} imply the following (see~\cite{MoOdOl:09} for more details):
\begin{theorem} \label{cor:bor}
Let $f,g : \R^n \to [0,1]$ be two measurable
functions on Gaussian space with $\E[f] = \mu$ and $\E[g] = \nu$.
Then for all $0 \leq \rho \leq 1$ we have
\[
\StabThrmax{\rho}{\mu,\nu} \leq \E[f U_{\rho} g] \leq \StabThrmax{\rho}{\mu,\nu}.
\]
\end{theorem}

We will need the following corollary.

\begin{corollary} \label{cor:bor_extended}
Let $X_1,\ldots,X_n,Y_1,\ldots,Y_m$ be jointly Gaussian such that $X_1,\ldots,X_n$ are independent, $Y_1,\ldots,Y_m$ are independent and
\[
\sup_{\|\alpha\|_2 = \|\beta\|_2 = 1} |\Cov[\sum_i \alpha_i X_i, \sum_i \beta_i Y_i]| \le\rho.
\]
Let $f_1 : \R^n \to [0,1], f_2 : \R^m \to [0,1]$
with $\E[f_j] = \mu_j$. Then
\[
\StabThrmin{\rho}{\mu_1,\mu_2}
\leq \E[f_1 f_2] \leq
\StabThrmax{\rho}{\mu_1,\mu_2}.
\]
\end{corollary}

\begin{proof}

Note that if $m > n$ then we may define $X_{n+1},\ldots,X_{m}$ to be Gaussian and
independent of all the other variables. This implies that without loss of generality we may assume that $m=n$.

We claim that without loss of generality the covariance matrix between $X$ and $Y$ is given by
\[
\Cov[X_i,Y_j] = 0,
\] for $i \neq j$ and
$|\Cov[X_i,Y_i]| \leq \rho$ for all $i$.
To see this, take $\| \alpha \|_2 = \| \beta \|_2 = 1$ which maximizes $\Cov[\sum \alpha_i X_i,\sum_i \beta_i Y_i]$ and define
$\tilde{X}_1 = \sum \alpha_i X_i, \tilde{Y}_1 = \sum \beta_i Y_i$ and $\tilde{X}_2,\ldots \tilde{X}_n$
to be orthonormal basis of the projection of the span of
$\{X_2,\ldots,X_n\}$
to the orthogonal complement of the space spanned by $\tilde{X}_1$ and similarly
$\tilde{Y}_2,\ldots,\tilde{Y}_n$. It is easy to check that
$\Cov[\tilde{X}_1,\tilde{Y}_i] = \Cov[\tilde{Y}_1,\tilde{X}_i] = 0$ for
$i > 1$. Repeating this process for $\tilde{X}_2,\ldots,\tilde{X}_n$ and $\tilde{Y}_2,\ldots,\tilde{Y}_n$ etc.
we obtain the desired covariance matrix.

Write $\Cov[X,Y] = \rho J$, where $J$ is diagonal with all entries in $[-1,1]$.
Then clearly
\[
\E[f_1(X_1,\ldots,X_n) f_2(Y_1,\ldots,Y_n)] = \E[f_1 U_{\rho} (U_J f_2)],
\]
where $U_{\rho}$ is the Ornstein-Uhlenbeck operator and $U_J$ is the operator
defined by
\[
(U_J f)(x_1,\ldots,x_n) =
\E[f(j_1 x_1 + \sqrt{1-J_{1,1}^2}\, y_1\,,\,\ldots,
         j_n x_n + \sqrt{1-J_{n,n}^2}\, y_n)],
\]
where $y$ is distributed according to the Gaussian measure.
Since $U_J$ is a Markov operator, we have $\E[U_J f_2] = \E[f_2]$ and
$0 \leq U_J f_2 \leq 1$.
Now applying Borell's result we obtain the desired conclusion.
\end{proof}

We note that in general there is no
closed form for $\StabThr{\rho}{\mu,\nu}$; however, for balanced functions we have
Sheppard's formula~\cite{Sheppard:99}: $\StabThr{\rho}{1/2,1/2} = \frac14+
\frac{1}{2\pi}\arcsin \rho$. Finally we record a fact to be used later.
\begin{proposition}
Let $((X_{i,j},Y_{i,j})_j)_{i=1}^n$
be jointly Gaussian, each distributed $N(0,1)$.
Suppose further that for each $i$:
\[
\sup_{\|\alpha\|_2 = \|\beta\|_2 = 1} |\Cov[\sum_j \alpha_j X_{i,j} \sum_i \beta_j Y_{i,j}]| \leq \rho,
\]
and that the $n$ collections $((X_{i,j},Y_{i,j})_j)_{i=1}^n$ are independent.
Then we have:
\[
\sup_{\|\alpha\|_2 = \|\beta\|_2 = 1} |\Cov[\sum_{i,j} \alpha_{i,j} X_{i,j} \sum_{i,j}
\beta_{i,j} Y_{i,j}]| \leq \rho,
\]
\end{proposition}

\begin{proof}
Using the fact that a linear combination of Gaussians is a Gaussian it suffices
to show that if $(X_i,Y_i)$ are independent Gaussian vectors, each satisfying
$|\Cov(X_i,Y_i)| \leq \rho, X_i, Y_i \sim N(0,1)$ and
$\| \alpha \|_2 = \| \beta \|_2 = 1$ then
\[
|\Cov[\sum_i \alpha_i X_i \sum_i \beta_i Y_i]| \leq \rho.
\]
This follows immediately from Cauchy-Schwarz:
\[
|\Cov[\sum_i \alpha_i X_i \sum_i \beta_i Y_i]| =
|\sum_i \alpha_i \beta_i \Cov[X_i,Y_i]| \leq \rho
\sum_i |\alpha_i \beta_i| \leq \rho.
\]

\end{proof}

\subsection{Asymptotics of $\StabThrmax{}{}$}
In some of the applications below we will need to estimate
$\StabThrmax{\rho_1,\ldots,\rho_{k-1}}{\mu_1,\ldots,\mu_k}$.
In particular, we will need the following estimate
\begin{lemma} \label{lem:rho_bound}
Let $0 < \mu < 1$ and $0 < \rho < 1$.
Define
\[
B_k(\rho,\mu) = \StabThrmax{\rho_1,\ldots,\rho_{k-1}}{\mu_1,\ldots,\mu_k},
\]
where $\rho_1 = \ldots = \rho_{k-1} = \rho$ and $\mu_1 = \ldots = \mu_k = \mu$.
Then as $k \to \infty$ we have
\begin{equation} \label{eq:Bbound}
B_k(\rho,\mu) \leq  k^{(\rho^2-1)/\rho^2 + o(1)}.
\end{equation}
\end{lemma}

\begin{proof}
Clearly, we have
\begin{equation} \label{eq:B_REC}
B_{i+1}(\rho,\mu) = \StabThrmax{\rho}{\mu,B_i}.
\end{equation}
The proof proceeds by deriving bounds on recursion~(\ref{eq:B_REC}).
This is a straightforward (but not very elegant) calculation with Gaussians.
Writing $B_i$ for $B_i(\rho,\mu)$, the main two steps in verifying~(\ref{eq:Bbound}) are to show that
\begin{itemize}
\item
The sequence
$B_j$ converges to $0$ as $j \to \infty$. This follows from the fact that
the functions $B \to \StabThr{\rho}{\mu,B}$ are easily seen to be
strictly decreasing and have no fixed points other than $B=0$ and $B=1$
when $0 < \rho < 1$.
\item
Using Gaussian estimates sketched below, we see that
for  $B_{j-1}$ sufficiently small it holds that
\begin{equation} \label{eq:finite_diff_B}
B_{j-1} - B_j \geq \frac{1}{2} B_{j-1}^{1 + \frac{\rho^2}{1-\rho^2} + o(1)}.
\end{equation}
This corresponds to $B_j$ of the form
\[
B_j \leq C j^{-\frac{1-\rho^2}{\rho^2} + o(1)}.
\]
More formally, it is easy to see that if $B_{j-1}$
is sufficiently small and satisfies
\[
B_{j-1} \leq C(j-1)^{-\alpha}
\]
and
\[
B_j \leq B_{j-1} (1-\frac{1}{2}B_{j-1}^{1/\alpha})
\]
for $\alpha > 0$
then
\[
B_{j} \leq C j^{-\alpha}.
\]
This follows using the fact that for small values of $\delta$
the maximum of the
function $x(1-x^{1/\alpha}/2)$ in the interval $[0,\delta]$
is obtained at $\delta$
and therefore:
\begin{eqnarray*}
B_{j} &\leq& B_{j-1}(1-\frac{1}{2}B_{j-1}^{1/\alpha}) \leq
C(j-1)^{-\alpha} \left( 1- \frac{1}{2}(C (j-1)^{-\alpha})^{1/\alpha} \right) \\ &=&
C(j-1)^{-\alpha} \left(1-\frac{1}{2}C^{1/\alpha} (j-1)^{-1} \right).
\end{eqnarray*}
In order that
\[
B_{j} \leq C j^{-\alpha},
\]
we need that
\[
\frac{1}{1 - \frac{1}{2} C^{1/\alpha} (j-1)^{-1}} \geq
\left( \frac{j}{j-1} \right)^{\alpha}
\]
which holds for large enough value of $C$.
\end{itemize}
In order to obtain~(\ref{eq:finite_diff_B}) for small values of $B_{j-1}$,
one uses the lemma stated below together with
the approximation
\[
\Phi(-t) = \exp(-(1+o_t(1)) \frac{t^2}{2}).
\]
which implies that for every fixed $\eps > 0$:
\[
\Phi(-t)
\Phi(-\frac{\rho(t+\eps)}{\sqrt{1-\rho^2}})  =
\Phi(-t)^{(1+o_t(1))(1+\frac{\rho^2}{1-\rho^2})}.
\]
\end{proof}

\begin{lemma}
Let $(X,Y)$ be a bi-variate normal with $\Var[X] = \Var[Y] = 1$ and
$\Cov[X,Y] = \rho$. Then for all $\eps > 0$ and $t > 0$ it holds that
\begin{eqnarray*}
\Phi(-t) - \StabThrmax{\rho}{1/2,\Phi(-t)} &=&
\P[X \leq -t] - \P[X \leq -t, Y \leq 0] = \P[X \leq -t, Y \geq 0] \\ &\geq&
\Phi(-t) (1-\exp(-t \eps - \eps^2/2))
\Phi(-\frac{\rho(t+\eps)}{\sqrt{1-\rho^2}})
\end{eqnarray*}
\end{lemma}

\begin{proof}
The equalities follow by the definitions.
For the inequality, we write
\begin{eqnarray*}
\P[X \leq -t, Y \geq 0] &\geq& \P[X \leq -t] \P[Y \geq 0, X \geq -t-\eps | X \leq -t] \\ &=&
\P[X \leq -t] \P[X \geq -t-\eps | X \leq -t] \P[Y \geq 0 | -t \geq X \geq -t-\eps].
\end{eqnarray*}
The bound in the lemma follows by bounding each of the three terms starting with
$\P[X \leq -t] = \Phi(-t)$.
Then note that
\[
\frac{\P[X \leq -t-\eps]}{\P[X \leq -t]} \leq
\exp(-t \eps - \eps^2/2),
\]
and therefore
\[
\P[X \geq -t-\eps | X \leq -t] \geq (1-\exp(-t \eps - \eps^2/2)).
\]
Finally, writing $Z$ for a $N(0,1)$ variable that is independent of
$X$ we obtain
\begin{eqnarray*}
\P[Y \geq 0 |-t \geq X \geq -t-\eps] \geq \P[Y \geq 0 | X = -t-\eps] &=&
\P[\rho (-t-\eps) + \sqrt{1-\rho^2}Z \geq 0] \\ &=&
\P[Z \geq \frac{\rho(t+\eps)}{\sqrt{1-\rho^2}}]  =
\Phi(-\frac{\rho(t+\eps)}{\sqrt{1-\rho^2}}),
\end{eqnarray*}
as needed.
\end{proof}

\section{Gaussian Bounds on Non-reversible Noise forms} \label{sec:conj}
In this section we prove the main results of the paper:
Theorem~\ref{thm:MIST_easy} and its relaxations Propositions~\ref{prop:MIST_relaxed_easy} and~\ref{prop:MIST_easy}.
As in previous work~~\cite{MoOdOl:05,MoOdOl:09,DiMoRe:06}, the proof idea is to use an invariance principle, in this case Theorem~\ref{thm:smooththeorem},
together with the Gaussian bounds of Section~\ref{sec:Gaussian}.

Since the invariance principle requires working either with low degree polynomials, or polynomials that have exponentially decaying weight,
an important step of the proof is the reduction to this case. This
reduction is proved in subsection~\ref{subsec:trunc}. It is based on the fact
that $\rho < 1$ and on the properties of correlated spaces and Efron-Stein
decompositions derived in Section~\ref{sec:correlated}.

The reduction, Theorem~\ref{thm:smooththeorem} and truncation arguments allow
to prove Theorem~\ref{thm:MIST_easy} for $k=2$ in subsection~\ref{subsec:maj_stablest} and for $k>2$ in subsection~\ref{subsec:multi_maj}.

The relaxed conditions on the influences for $k=2$ and for $r$-wise
independent distributions are derived in subsection~\ref{subsec:relaxed_inf}
using a ``two-threshold'' technique. A related technique has been used before
in~\cite{DinurSafra:02,DiMoRe:06}. However, the variant presented here is more
elegant, gives more explicit dependency on the influences and allows to exploit
$s$-wise independence.

Finally, using a recursive argument
we derive in subsection~\ref{subsec:reg}
Proposition~\ref{prop:MIST_relaxed_easy}.
We don't know of any previous application of this idea in the context of the study
of influences.

\subsection{Noise forms are determined by low degree expansion} \label{subsec:trunc}
In order to use the invariance principle, it is crucial to apply it to
multi-linear polynomials that are either of low degree or well approximated by
their low degree part. Here we show that noise stability quantities do not
change by much if one replaces a function by a slight smoothing of it.
For the following statement recall Definition~\ref{def:rho} for the definition of $\rho$ and (\ref{eq:bonami_beckner1}),(\ref{eq:bonami_beckner2}) and~(\ref{eq:bonami_beckner3}) for the definition of the Bonami-Beckner operator $T_{1-\gamma}$.

\begin{lemma} \label{lem:smooth}
Let $\Omega_1,\ldots,\Omega_n,\Lambda_1,\ldots,\Lambda_n$ be a collection
of finite probability spaces.
Let $\CalX, \CalY$ be two ensembles such the collections of random variables
such that $(\CalX_i,\CalY_i)$ are independent and $\CalX_i$ is a basis for the
functions  in $\Omega_i$, $\CalY_i$ is a basis for the functions in $\Lambda_i$.
Suppose further that for all $i$ it holds that
$\rho(\Omega_i,\Lambda_i) \leq \rho$.

Let $P$ and $Q$ be two multi-linear polynomials. Let $\eps > 0$ and $\gamma$ be chosen sufficiently close to $0$ so that
\[
\gamma \leq 1 - (1-\eps)^{\log \rho / (\log \eps + \log \rho)}.
\]
Then:
\[
| \E[P(\CalX) Q(\CalY)] - \E[T_{1-\gamma} (\CalX) T_{1-\gamma}Q(\CalY)] | \leq
2 \eps \Var[P] \Var[Q].
\]
In particular, there exists an absolute constant $C$ such that it suffices to take
\[
\gamma = C \frac{(1-\rho)\eps}{\log (1/\eps)}.
\]
\end{lemma}

\begin{proof}
Without loss of generality if suffices to assume that
$\Var[P] = \Var[Q] = 1$ and show that
\[
| \E[P(\CalX) Q(\CalY)] - \E[P (\CalX) (T_{1-\gamma} Q)(\CalY)] | =
| \E[P(\CalX) \left((I - T_{1-\gamma}) Q \right) (\CalY)] | \leq \eps.
\]
Let $T$ be the Markov operator defined by
$Tg(x) = \E[g(Y) | X = x]$, where $(X,Y)$
are distributed according to $(\CalX,\CalY)$.
In order to prove the lemma it suffices to show that
\begin{equation} \label{eq:trunc_lemma}
| \E[P(\CalX) \left(T(I - T_{1-\gamma}) Q \right) (\CalX)] | \leq \eps.
\end{equation}

Write $P$ and $Q$ in terms of their Efron-Stein decomposition, that is,
\[
P = \sum_{S \subset [n]} P_S, \quad
Q = \sum_{S \subset [n]} P_Q.
\]
It is easy to see that
\[
(I-T_{1-\gamma}) Q_S = (1-(1-\gamma)^{|S|}) Q_S,
\]
and propositions~\ref{prop:TEfronStein1} and~\ref{prop:TEfronStein2} imply
that
\[
\| T Q_S \|_2 \leq \rho^{|S|} \| Q_S \|_2,
\]
and that  $TQ_S$ is orthogonal to $P_{S'}$ for $S' \neq S$.
Writing $T' = T(I - T_{1-\gamma})$ we conclude that
\[
\| T' Q_S \|_2 \leq \min(\rho^{|S|}, 1-(1-\gamma)^{|S|}) \| Q_S \|_2
\leq \eps \| Q_S \|_2,
\]
and that $T' Q_S$ is orthogonal to $P_S'$ when $S' \neq S$.

By Cauchy-Schwarz we get that
\[
| \E[P(\CalX) (T'Q) (\CalX)] | =
|\sum_{S \neq \emptyset} \E[P_S T' Q_S]|
\leq \sqrt{\Var[P]} \sqrt{ \sum_{S \neq \emptyset} \| T' Q_S \|_2^2}
\leq \eps \sqrt{\Var[P]} \sqrt{\Var[Q]},
\]
as needed.

\end{proof}

Similarly we have:
\begin{lemma} \label{lem:smooth_general}
Let $(\Omega_1^{(j)},\ldots,\Omega_n^{(j)})_{j=1}^k$ be $k$ collections
of finite probability spaces.
Let $( (\CalX_i^j)_{i=1}^n : j=1,\ldots,k)$ be
$k$ ensembles of collections of random variables
such that $((\CalX_i^j)_{j=1}^k)_{i=1}^n$
are independent and $\CalX_i^j$ is a basis for the
functions  in $\Omega_i^{(j)}$.
Suppose further that for all $i$ it holds that
$\rho(\Omega_i^{(j)} : 1 \leq j \leq k) \leq \rho$.

Let $P_1,\ldots,P_k$ be $k$ multi-linear polynomials.
Let $\gamma$ be chosen sufficiently close to $0$ so that
\[
\gamma \leq 1 - (1-\eps)^{\log \rho / ( \log \eps + \log \rho)}.
\]
Then:
\[
| \E[\prod_{j=1}^k P_j(\CalX^j)] -
  \E[\prod_{j=1}^k T_{1-\gamma} P_j (\CalX^j)]| \leq \eps
\sum_{j=1}^k \sqrt{\Var[P_j]}
\sqrt{\Var[\prod_{j' < j} T_{1-\gamma} P_j \prod_{j' > j} P_j]}.
\]
In particular, there exists an absolute constant $C$ such that it suffices to take
\[
\gamma = C \frac{(1-\rho)\eps}{\log (1/\eps)}.
\]
\end{lemma}

\begin{proof}
The proof follows the proof of the previous lemma.
\end{proof}

\subsection{Bilinear Gaussian Bounds} \label{subsec:maj_stablest}
In this section we prove the bilinear stability bound.
We repeat the statement of Theorem~\ref{thm:MIST_easy}
with more explicit dependency on the influences.
\begin{theorem} \label{thm:MIST}
Let $(\Omega_i^{(1)} \times \Omega_i^{(2)}, \P_i)$ be a sequence of correlated
spaces such that for each $i$ the minimum $\P_i$ probability of any atom in
is at least $\alpha \leq 1/2$ and such that $\rho(\Omega_i^{(1)},\Omega_i^{(2)} ; \P_i) \leq \rho$ for all $i$.
Then for every $\eps > 0$ there exists a
$\tau = \tau(\eps) < 1/2 $
that if
$f: \prod_{i=1}^n \Omega_i^{(1)} \to [0,1]$ and
$g: \prod_{i=1}^n \Omega_i^{(2)} \to [0,1]$ satisfy
\begin{equation} \label{eq:inf_trunc_two}
\max(\Inf_i^{\leq \log(1/\tau) \log(1/\alpha)}(f),  \Inf_i^{\leq \log(1/\tau) \log(1/\alpha)}(g) ) \leq \tau
\end{equation}
for all $i$ (See
Definition~\ref{def:low-degree-influence} for the definition of
low-degree influence) then
\begin{equation} \label{eq:corr_bd3}
\StabThrmin{\rho}{\E[f],\E[g]} - \eps
\leq \E[f g] \leq
\StabThr{\rho}{\E[f],\E[g]} + \eps.
\end{equation}
Moreover, there exists an absolute constant $C$ such that one may take
\begin{equation} \label{eq:tau_bd3}
\tau = \eps^{C \frac{\log(1/\alpha) \log(1/\eps)}{\eps (1-\rho)}}
\end{equation}
\end{theorem}

\begin{proof}
Write $\mu = \E[f]$ and $\nu = \E[g]$ and $K = \log(1/\alpha)$.
As discussed in Section~\ref{sec:setup}, let $\tilde{\CalX}$ be the
sequence of ensembles such that $\tilde{\CalX}_i$ spans the functions on
$\Omega_i = \Omega_i^{(1)} \times \Omega_i^{(2)}$, $
\CalX_i$ spans functions on $\Omega_i^{(1)}$ and $\CalY_i$ spans
the functions on $\Omega_i^{(2)}$. We now express $f$ and $g$ as
multilinear polynomials $P$ and $Q$ of $\CalX$ and $\CalY$.
Let $\gamma > 0$ be chosen so that
\[
| \E[P(\CalX) Q(\CalY)] - \E[T_{1-\gamma} (\CalX) T_{1-\gamma}Q(\CalY)] | \leq
  \eps/2.
\]
Note that by Lemma~\ref{lem:smooth} it follows that we may take
$\gamma = \Theta(\eps (1-\rho) / \log (1/\eps))$.
Thus it suffices to prove the bound stated in the theorem for
$T_{1-\gamma}P(\CalX)$ and $T_{1-\gamma} Q(\CalY)$.

We use the invariance principle under hypothesis $\hthree$.
Let $\CalG$ and $\CalH$ be two Gaussian ensembles such that for all $i$ the
covariance matrix of $\CalH_i$ and $\CalG_i$ is identical to the covariance
matrix of $\CalX_i$ and $\CalY_i$ and such that $(\CalG_i,\CalH_i)$
are independent. Clearly:
\[
\E[T_{1-\gamma}P (\CalX) \, T_{1-\gamma} Q(\CalY)] =
\E[T_{1-\gamma}P (\CalG) \, T_{1-\gamma} Q(\CalH)].
\]

Since $(P(\CalX), Q(\CalY))$ takes values in $[0,1]^2$
the same is true for
\[
(\tilde{P},\tilde{Q}) = ((T_{1-\gamma} P)(\CalX),T_{1-\gamma} Q(\CalY)).
\]
In other words, $\E[\trunc(\tilde{P},\tilde{Q})] = 0$,
where $\trunc$ is the function in~(\ref{eq:def_trunc}).
Writing
\[
(\bar{P},\bar{Q}) = ((T_{1-\gamma} P)(\CalG),T_{1-\gamma}Q(\CalH)),
\]
we conclude from Theorem~\ref{thm:smooththeorem} that
$\E[\trunc(\bar{P},\bar{Q}))] \leq \tau^{\Omega(\gamma/K)}$.
That is, $\|(\bar{P},\bar{Q}) - (P',Q')\|_2^2 \leq
\tau^{\Omega(\gamma/K)}$, where $P''$ is the function of $\bar{P}$ defined by
\[
P' = \left\{ \begin{array}{rl}
                0    & \text{if $\bar{P} \leq 0$,} \\
                \bar{P} & \text{if $\bar{P} \in [0,1]$,} \\
                1    & \text{if $\bar{P} \geq 1$.}
                \end{array}
                \right. ,
\]
and $Q'$ is defined similarly.
Now using Cauchy-Schwarz it is easy to see that
\[
|\E[\bar{P} \bar{Q}] - \E[P' Q']| \leq \tau^{\Omega(\gamma/K)}.
\]
Write $\mu' = \E[P']$ and $\nu' = \E[Q']$. Then using the Gaussian Corollary~\ref{cor:bor} we obtain that
\[
\E[P' Q'] \leq \StabThr{\rho}{\mu',\nu'}.
\]

From Cauchy-Schwarz it follows that
$|\mu-\mu'| \leq \tau^{\Omega(\gamma/K)}$ and similarly for $\nu,\nu'$.
It is immediate to check that
\[
|\StabThr{\rho}{\mu,\nu} -
\StabThr{\rho'}{\mu',\nu'}| \leq |\mu-\mu'| + |\nu-\nu'| \leq
\tau^{\Omega(\gamma/K)}.
\]
Thus we have
\[
\E[PQ] \leq \StabThr{\rho}{\mu,\nu} + \tau^{\Omega(\gamma/K)} + \eps/2.
\]
Taking $\tau$ as in~(\ref{eq:tau_bd3}) yields
\[
\tau^{\Omega(\gamma/K)} \leq \eps/2,
\]
and thus we obtain the upper bound in~(\ref{eq:corr_bd3}). The proof of the
lower bound is identical.
\end{proof}

The following proposition completes the proof of~(\ref{eq:corr_bd}).
\begin{proposition} \label{prop:MIST}
Let $(\prod_{j=1}^k \Omega_i^{(j)},\P_i), 1 \leq i \leq n$
be a sequence of correlated spaces such that for each $i$
the minimum $\P_i$ probability of any atom in
$\Omega_i$ is at least $\alpha \leq 1/2$ and
such that $\rho(\Omega^{(1)}_i,\ldots,\Omega^{(k)}_i) \leq \rho < 1$ for all $i$.

Then for every $\eps > 0$ there exists a $\tau > 0$
such that if $f_j: \prod_{i=1}^n \Omega^{(j)}_i \to [0,1]$ for $1 \leq j \leq k$ satisfy
\[
\max(\Inf_i^{\leq \log(1/\tau)/\log(1/\alpha)}(f_j)) \leq \tau
\]
for all $i$ and $j$ (See
Definition~\ref{def:low-degree-influence} for the definition of
low-degree influence) then it holds that
\[
\StabThrmin{\underline{\rho}}{\E[f_1],\ldots,\E[f_k]} - \eps
\leq \E[\prod_{i=1}^k f_i] \leq
\StabThrmax{\underline{\rho}}{\E[f_1],\ldots,\E[f_k]} + \eps
\]
There exists an absolute constant $C$ such that one may take
\begin{equation} \label{eq:tau_bd}
\tau = \eps^{C \frac{\log(1/\alpha) \log(1/\eps)}{\eps (1-\rho)}}
\end{equation}
\end{proposition}

The proof uses the following lemma, see e.g.~\cite{SamorodnitskyTrevisan:06}.
\begin{lemma}
Let $f_1,\ldots,f_k : \Omega^n \to [0,1]$. Then for all $j$:
\[
\Inf_j(\prod_{i=1}^k f_i) \leq k \sum_{i=1}^k \Inf_j(f_i).
\]
\end{lemma}

\begin{proof}
The proof is based on the fact that
\[
\Var[\prod_{i=1}^k f_i] = \frac{1}{2}\E[(\prod_{i=1}^k f_i (X) - \prod_{i=1}^k f_i(Y))^2],
\]
where $X$ and $Y$ are independent. Now
\[
\E[(\prod_{i=1}^k f_i (X) - \prod_{i=1}^k f_i(Y))^2] \leq
\E[(\sum_{i=1}^k |f_i (X) - f_i(Y)|)^2] \leq
k \sum_{i=1}^k \E[(f_i(X)-f_i(Y))^2] = 2k \sum_{i=1}^k \Inf_i(f)
\]
which gives the desired result.
\end{proof}

\begin{proof}[Of Proposition~\ref{prop:MIST}]
The proof follows by applying Theorem~\ref{thm:MIST}
iteratively for the functions
$f_1, f_1 f_2, f_1 f_2 f_3, \ldots$ and using the previous lemma and the fact that $\StabThrmax{\rho}{}$ and $\StabThrmin{\rho}{}$ have Lipchitz constant $1$ in each coordinate.
\end{proof}



\subsection{Multi-linear bounds} \label{subsec:multi_maj}
Next we prove~(\ref{eq:ind_bd}).

\begin{theorem} \label{thm:ARIT_STAB}
Let $(\prod_{j=1}^k \Omega_j^{(i)},\P_i), 1 \leq i \leq n$
be a sequence of correlated spaces such that for each $i$
the minimum $\P_i$ probability of any atom in
$\Omega_i$ is at least $\alpha \leq 1/2$ and
such that $\rho(\Omega^{(1)}_i,\ldots,\Omega^{(k)}_i) \leq \rho < 1$ for all $i$.
Suppose furthermore that for all $i$ and $j \neq j'$, it holds that
$\rho(\Omega^{(j)}_i,\Omega^{(j')}_i) = 0$, i.e.,
$(\prod_{j=1}^k \Omega^{(j)}_i, \P_i)$ is pairwise independent.

Then for every $\eps > 0$ there exists a
$\tau = \tau(\eps) < 1/2 $
such that if $f_j: \prod_{i=1}^n \Omega^{(j)}_i \to [0,1]$ for $1 \leq j \leq k$ satisfy
\begin{equation} \label{eq:arit_stab_inf}
\max(\Inf_i^{\leq \log(1/\tau)/\log(1/\alpha)}(f_j)) \leq \tau
\end{equation}
then
\begin{equation} \label{eq:corr_bd2}
|\E[\prod_{i=1}^k f_i] - \prod_{i=1}^k \E[f_i]| \leq \eps.
\end{equation}
There exists an absolute constant $C$ such that one may take
\begin{equation} \label{eq:tau_bd2}
\tau = \eps^{C \frac{K \log(1/\eps)}{\eps (1-\rho)}}
\end{equation}

\end{theorem}

\begin{proof}
Note that for all $i$ and all $\gamma$ we have that the functions
$f_i$ and $T_{1-\gamma} f_i$ are $[0,1]$ valued functions. Therefore, as
in the previous proof we obtain by Lemma~\ref{lem:smooth_general} that
\[
|\E[\prod_{i=1}^k f_i] - \E[\prod_{i=1}^k T_{1-\gamma} f_i]| \leq \eps/4.
\]
for
$\gamma = \Omega(\eps (1-\rho) / \log (1/\eps))$.
Thus it suffices to prove the bound stated in the theorem for the functions
$T_{1-\gamma} f_i$.

We now use the invariance principle. Recall that the $f_i$ may be written
as a multi-linear polynomial $P_i$ of an ensemble $\CalX_i$.
Let $\CalG_i, 1 \leq i \leq k$
denote Gaussian ensembles with the same covariances as
$\CalX_i, 1 \leq i \leq k$ and  let $(g_1,\ldots,g_k)$ be the multi-linear
polynomials $(P_1,\ldots,P_k)$ applied to $(\CalX_1,\ldots,\CalX_k)$.
Let $h_i = f_{[0,1]}(T_{1-\gamma} g_i)$.
By the invariance principle
Corollary~\ref{cor:smooththeorem}, we have:
\[
|\E[\prod_{i=1}^k T_{1-\gamma} f_i] - \E[\prod_{i=1}^k h_i]| \leq \eps/4.
\]
Note that $h_i$ are functions of ensembles of Gaussian random variables such that each pair of ensembles is independent. Therefore, the $h_i$'s are
independent which implies in turn that
\[
\E[\prod_{i=1}^k h_i] = \prod_{i=1}^k \E[h_i].
\]
Corollary~\ref{cor:smooththeorem} also implies that
\[
|\E[h_i] - \E[f_i]| = |\E[h_i] - \E[T_{1-\gamma} f_i]| \leq \eps/4k,
\]
so
\[
|\prod_{i=1}^k \E[h_i] - \prod_{i=1}^k \E[f_i]| \leq \eps/4,
\]
which concludes the proof.
\end{proof}

\subsection{Relaxed influence conditions} \label{subsec:relaxed_inf}
In this subsection we will relax the conditions imposed on the influence,
i.e. Proposition~\ref{prop:MIST_relaxed_easy}.
In particular we will show that in Theorem~\ref{thm:MIST_easy} and $k=2$
it suffices to assume that for each coordinate at least one of the
functions has low influence. Similarly for $k>2$ and $s$-wise independent
distributions it suffices to have that in each coordinate at most
$s$ of the functions have large influence.

\begin{lemma} \label{lem:average}
Assume $\mu$ is an $r$-wise independent distribution on $\Omega^k$.
Let $f_1,\ldots,f_k : \Omega^n \to [0,1]$. Let $S \subset [n]$ be a set
of coordinates such that for each $i \in S$ at most $r$ of the functions
$f_j$ have $\Inf_i(f_j) > \eps$. Define
\[
g_i(x) = \E[f_i(Y) | Y_{[n] \setminus S} = x_{[n] \setminus S}].
\]
Then the functions $g_i$ do not depend on the coordinates in $S$ and
\[
|\E[\prod_{i=1}^k f_i] - \E[\prod_{i=1}^k g_i]| \leq k |S| \sqrt{\eps}.
\]
\end{lemma}

\begin{proof}
Recall that averaging over a subset of the variables preserves expected value.
It also maintains the property of taking values in $[0,1]$
and decreases influences. Thus it suffices to prove the claim for the case
where $|S| = 1$. The general case then follows by induction.

So assume without loss of generality that $S=\{1\}$ consists of the first
coordinate only and that $f_{r+1},\ldots,f_k$ all have $\Inf_1(f_j) \leq \eps$, so
that $\E[(f_j-g_j)^2] \leq \eps$ for $j > r$. Then by Cauchy-Schwarz we have
$\E[|f_j-g_j|] \leq \sqrt{\eps}$ for $j > r$ and
using the fact that the functions are bounded in $[0,1]$ we obtain
\begin{equation} \label{eq:eps_r_wise}
|\E[\prod_{i=1}^k f_i] - \E[\prod_{i=1}^r f_i \prod_{i=r+1}^k g_i]| \leq
\E[|\prod_{i=r+1}^k f_i - \prod_{i=r+1}^k g_i |] \leq
\sum_{i=r+1}^k \E[|f_i - g_i|] \leq k \sqrt{\eps}.
\end{equation}
Let us write $\E_1$ for the expected value with respect to the first variable.
Recalling that the $g_i$ do not depend on the first variable and
that the $f_i$ are $r$-wise independent we obtain that
\[
\E_1[\prod_{i=1}^r f_i \prod_{i=r+1}^k g_i]] =
\prod_{i=r+1}^k g_i
\E_1[\prod_{i=1}^r f_i] =
\prod_{i=r+1}^k g_i
\prod_{i=1}^r \E_1[f_i] =
\prod_{i=1}^r g_i \prod_{i=r+1}^k g_i = \prod_{i=1}^k g_i.
\]
This implies that
\begin{equation} \label{eq:prod_r_wise}
\E[\prod_{i=1}^r f_i \prod_{i=r+1}^k g_i]] = \E[\prod_{i=1}^k g_i],
\end{equation}
and the proof follows from~(\ref{eq:eps_r_wise}) and(\ref{eq:prod_r_wise}).
\end{proof}

We now prove that condition~(\ref{eq:inf_bd_condition_relaxed1}) suffices
instead of~(\ref{eq:inf_bd_condition}).\\

\begin{lemma}
Theorem~\ref{thm:MIST} holds with the condition
\begin{equation} \label{eq:inf_trunc_two_relaxed}
\max_{i}(\min(\Inf_i^{\leq \log^2(1/\tau)/\log(1/\alpha)}(f_1),\Inf_i^{\leq \log^2(1/\tau)/\log(1/\alpha)} (f_2)))
\leq \tau
\end{equation}
instead of~(\ref{eq:inf_trunc_two}).
\end{lemma}

\begin{proof}
Let
\[
\gamma = C_1 \frac{(1-\rho) \eps}{\log(1/\eps)},
\]
\[
\tau = \eps^{C_2 \frac{\log(1/\alpha) \log(1/\eps)}{(1-\rho)\eps}},
\]
and
\[
R = \frac{\log(1/\tau)}{\log(1/\alpha)}.
\]
From the proof of Theorem~\ref{thm:MIST} it follows that the constant $C_1$ and $C_2$ may be chosen
so that
\[
(1-\gamma)^R \leq \eps/4,
\]
\[
|\E[f g] - \E[T_{1-\gamma} f T_{1-\gamma} g]| \leq \eps/4.
\]
Moreover the conclusion of Theorem~\ref{thm:MIST} holds with error at most $\eps/4$
for any pair of functions $\tilde{f}$ and $\tilde{g}$
if {\em all} influences of $\tilde{f}$ and $\tilde{g}$ satisfy $\Inf_i^{\leq R} \leq \tau$.
Let
\[
R' = \frac{\log^2(1/\tau)}{\log(1/\alpha)}
\]
and choose $C_3$ large enough so that
\[
\tau' = \eps^{C_3 \frac{\log(1/\alpha) \log(1/\eps)}{(1-\rho)\eps}},
\]
satisfy
\[
\frac{2 R}{\tau} \sqrt{\tau' + (1-\gamma)^{2 R'}} \leq \frac{\eps}{4}.
\]

Assume that $f$ and $g$ satisfy
\[
\max_{i}(\min(\Inf_i^{\leq R'}(f),\Inf_i^{\leq R'}(g))) \leq \tau'.
\]
We will show that the statement of the theorem holds for $f$ and $g$.
For this let
\[
S_f = \{ i : \Inf_i^{\leq R}(f) > \tau \}, \quad
S_g = \{ i : \Inf_i^{\leq R}(g) > \tau \}.
\]
Let $S = S_f \cup S_g$.
Since $R' \geq R$ and $\tau' \leq \tau$, the sets $S_f$ and $S_g$ are disjoint.
Moreover, both $S_f$ and $S_g$ are of size at most $\frac{R}{\tau}$.
Also, if $i \in S$ and $\Inf_i^{\leq R}(f) > \tau$ then
$\Inf_i^{\leq R'}(g) < \tau'$ and therefore $\Inf_i(g) \leq \tau' + (1-\gamma)^{2 R'}$.
In other words, for all $i \in S$ we have
$\min(\Inf_i(g),\Inf_i(f)) \leq \tau' + (1-\gamma)^{2 R'}$.

Letting
$S = S_f \cup S_g$ and applying Lemma~\ref{lem:average} with
\[
g'(y') = \E[g(Y) | Y_{[n] \setminus S} = y'_{[n] \setminus S}], \quad
f'(x') = \E[f(X) | X_{[n] \setminus S} = x'_{[n] \setminus S}],
\]
we obtain that
\begin{equation} \label{eq:trunc_double_coef}
|\E[f(x) g(y)] - \E[f'(x) g'(y)]| \leq
\frac{2 R}{\tau} \sqrt{\tau' + (1-\gamma)^{2 R'}}  \leq \frac{\eps}{4}.
\end{equation}
Note that the functions $f'$ and $g'$ satisfy that
$\max(\Inf_i^{\leq R}(f'),\Inf_i^{\leq R}(g')) \leq \tau$ for all $i$.
This implies that the results of Theorem~\ref{thm:MIST}
hold for $f'$ and $g'$. This together with~(\ref{eq:trunc_double_coef}) implies
the desired result.
\end{proof}

The proof of the relaxed condition on influences in Theorem~\ref{thm:ARIT_STAB}
is similar.
\begin{lemma}
Assume the setup of Theorem~\ref{thm:ARIT_STAB} where
$(\prod_{j=1}^k \Omega_i^{(j)}, \P_i)$ is $s$-wise independent for all $i$.
Then the conclusion of the theorem holds when the following condition:
\begin{equation} \label{eq:inf_trunc_two_relaxed2}
\forall i, | \{ j : \Inf_i^{\leq \log^2(1/\tau)/\log(1/\alpha)} (f_j) > \tau \} | \leq s
\end{equation}
replaces condition~(\ref{eq:arit_stab_inf}).
\end{lemma}

\begin{proof}
Again we start by looking at $T_{1-\gamma} f_i$
where $\gamma = \Omega(\eps (1-\rho) / \log (1/\eps))$.
We let $\tau'$ and $R'$ be chosen so that
\[
\frac{k R}{\tau} \sqrt{\tau' + (1-\gamma)^{2 R'}} \leq \frac{\eps}{4}.
\]
The set $S$ will consist of all coordinates $j$ where at least one
of the functions $f_i$ has $\Inf_j^{\leq R}(f_i) > \tau$. The rest of the proof is
similar.
\end{proof}

\subsection{A Recursive Argument} \label{subsec:reg}
Here we show how Proposition~\ref{prop:MIST_easy} follows from
Theorem~\ref{thm:MIST_easy}. The proof uses the following lemma.

\begin{lemma} \label{lem:reg1}
Let $(\Omega_1, \mu_1), \ldots, (\Omega_n,\mu_n)$
be finite probability spaces such that for all
$i$ the minimum probability of any atom in $\Omega_i$ is at least
$\alpha \leq 1/2$.  Let $f : \prod_{i=1}^n \Omega_i \to [0,1]$ and suppose
that $\Inf_i(f) \geq \tau$. Let $S = \{i\}$. Then
\[
\E[\overline{f}^S]  \geq \E[f] + \alpha \tau , \quad
\E[\underline{f}^S] \leq \E[f] - \alpha \tau .
\]
\end{lemma}

\begin{proof}
Note that if $g: \Omega_i \to [0,1]$ satisfies $\Var[g] \geq a$
then
\[
(\max_x g(x) - \E[g])^2 \geq a \alpha, \quad
(\min_x g(x) - \E[g])^2 \geq a \alpha.
\]
Therefore:
\[
\E[\overline{f}^S-f] \geq
\E[(\overline{f}^S-f)^2] \geq \alpha \Inf_i(f) = \alpha \tau,
\]
which implies the first inequality. The second inequality is proved
similarly.
\end{proof}



We now prove Proposition~\ref{prop:MIST_easy}.

\begin{proof}
The set $T$ is defined recursively as follows.
We let $T_o = \emptyset$.  $\underline{a}_t = \overline{a}_t = 0$.
Then we repeat the following: If all of the functions
\[
\overline{f}_1^{T_t},\ldots,\overline{f}_k^{T_t},
\underline{f}_1^{T_t},\ldots,\underline{f}_k^{T_t}
\]
have influences lower than
$\tau$, then we halt and let $T = T_t$. Otherwise, there exists at least one
$i$ and one $j$ such that either $\Inf_i(\overline{f}_j^{T_t}) \geq \tau$
or $\Inf_i(\underline{f}_j^{T_t}) \geq \tau$.
We then
let $T_{t+1} = T_t \cup \{ i \}$.  In the first case we let
$\overline{a}_{t+1} = \overline{a}_t + 1, \underline{a}_{t+1} = \underline{a}_t$.
In the second case we let $\overline{a}_{t+1} = \overline{a}_t, \underline{a}_{t+1} = \underline{a}_t + 1$.

Note that by Lemma~\ref{lem:reg1},
this process must terminate within
\[
\frac{2k}{\alpha \tau}
\]
steps since
\[
k \geq \E[\sum_{j=1}^k \overline{f}_j^{T_t}] \geq \alpha \tau \overline{a}_t, \quad
0 \leq \E[\sum_{j=1}^k \underline{f}_j^{T_t}] \leq k-\alpha \tau \underline{a}_t.
\]
and $\overline{a}_t + \underline{a}_t \geq t$.
\end{proof}

\section{Applications to Social Choice} \label{sec:social}
In this section we apply Theorem~\ref{subsec:maj_stablest}
to the two social choice models.

\subsection{$\rho$ for samples of votes}
In the first social choice example we consider Example~\ref{ex:sample}.
The correlated probability spaces are the ones given by
\[
\Omega_V = \{ \{x=1\}, \{x=-1\}\},
\]
representing the intended vote and
\[
\Omega_S = \{ \{(x=1,y=1)\}, \{(x=-1,y=1)\}, \{y=0\} \}
\]
representing the sampled status.

In order to calculate $\rho(\Omega_1,\Omega_2)$ it suffices by
lemma~\ref{lem:rho_cond_exp} to calculate $\sqrt{E[(Tf)^2]}$
where $f(x,y) = x$ is the (only) $\Omega_V$ measurable with $\E[f] = 0$ and
$\E[f^2] = 1$.
We see that $Tf(x,y) = 0$ if $y = 0$ and $Tf(x,y) = x$ when $y \neq 0$.
Therefore
\[
\sqrt{\E[(Tf)^2]} = \rho^{1/2}.
\]
\begin{lemma} \label{lem:Omega_SV}
\[
\rho(\Omega_V,\Omega_S) = \rho^{1/2}.
\]
\end{lemma}

\subsection{Predictability of Binary Vote}
Here we prove Theorem~\ref{thm:MIMP}.\\

\begin{proof}
The proof follows directly from Proposition~\ref{prop:MIST_relaxed_easy} and
Lemma~\ref{lem:Omega_SV} as
\[
\StabThr{\sqrt{\rho}}{1/2,1/2} = \frac14+
\frac{1}{2\pi}\arcsin \sqrt{\rho}.
\]
\end{proof}

\subsection{$\rho$ in Condorcet voting}
In the context of Condorcet voting, $\Omega$ is given
by $S_k$, the group of permutations on $k$ elements and $\mu$ is the
uniform measure. We write ${[k] \choose 2}$ for the collection of subsets of $[k]$ of size $2$.
\begin{definition}
Let $Q \subset {[k] \choose 2}$ and define
$R_Q : \Omega \to \{0,1\}^Q$ be letting $(R_Q(\sigma))_{i < j} = 1$ if
$\sigma(i) < \sigma(j)$ and $(R_Q(\sigma))_{i < j} = 0$ if
$\sigma(i) > \sigma(j)$. $R_Q$ summarizes the pairwise relations
in the permutation $\sigma$ for pairs in $Q$.

Given a subset $Q \subset {[k] \choose 2}$ we define
\[
\Omega_Q = \{ \{ \sigma : R(\sigma) = x\} : x \in \{0,1\}^Q \}
\]
Thus $\Omega_Q$ is the coarsening of $\Omega$ summarizing the information
about pairwise relations in $Q$.

\end{definition}

We will mostly be interested in $\rho(\Omega_Q,\Omega_{i<j})$ where
$(i < j) \notin Q$.
\begin{lemma} \label{lem:rho_max_vote}
Suppose $Q = \{(1 > 2),(1 > 3),\ldots,(1 > r)\}$ then
\[
\rho(\Omega_Q,\Omega_{1 > (1+r)}) = \sqrt{\frac{r-1}{3(r+1)}} \leq
\frac{1}{\sqrt{3}}
\]
\end{lemma}

\begin{proof}
We use lemma~\ref{lem:rho_cond_exp} again.
The space $\Omega_{1 > (r+1)}$ has a unique function with $\E[f] = 0$ and
$\E[f^2] = 1$. This is the function that satisfies $f(\sigma) = 1$ when
$\sigma(1) > \sigma(r+1)$ and $f(\sigma) = -1$ if $\sigma(1) < \sigma(r+1)$.

The conditional probability that $\sigma(1) > \sigma(r+1)$ given that $s$ of
the inequalities $\sigma(1) > \sigma(2),\ldots,\sigma(1) > \sigma(r)$ hold is
\[
\frac{s+1}{r+1}.
\]
Therefore the conditional expectation of $f$ under this conditioning is:
\[
\frac{(s+1)-(r-s)}{r+1} = \frac{2s+1-r}{r+1}.
\]
Noting that the number
of inequalities satisfied is uniform in the range $\{0,\ldots,r-1\}$
we see that
\begin{eqnarray*}
\E[(Tf)^2] &=& \frac{1}{r (r+1)^2} \sum_{s=0}^{r-1} (2s+1-r)^2 =
\frac{1}{r (r+1)^2}
\left( 4 \sum_{s=0}^{r-1} s^2 - 4(r-1) \sum_{s=0}^{r-1} s + r (r-1)^2 \right) \\ &=&
\frac{1}{r (r+1)^2}
\left( 4 \frac{2(r-1)^3 + 3(r-1)^2 + (r-1)}{6}
       - 4(r-1) \frac{(r-1)^2 + (r-1)}{2} + r (r-1)^2 \right) \\ &=&
\frac{r-1}{3(r+1)}.
\end{eqnarray*}
Therefore:
\[
\rho = \sqrt{\frac{r-1}{3(r+1)}}.
\]
\end{proof}

\subsection{Condorcet Paradox}
We now prove Theorem~\ref{thm:mibfc} dealing with Condorcet paradoxes. \\

\begin{proof}

We wish to bound
the asymptotic probability in terms of the number of candidates $k$
for the probability that there is a unique maximum in Condorcet aggregation.
 Clearly this probability equals $k$ times the probability that candidate
$1$ is the maximum.

Recall that the votes are aggregated as follows.
Let $f : \bits^n \to \bits$ be an anti-symmetric be a function with $\E[f] = 0$ and low influences. Let $\sigma_1,\ldots,\sigma_n \in S_k$
denote the $n$ rankings of the individual voters. Recall that we denote
$x^{a > b}(i) = 1$ if $\sigma_i(a) > \sigma_i(b)$ and $x^{a < b}(i) = -1$ if
$\sigma_i(a) < \sigma_i(b)$. Note that $x^{b > a } = -x^{a > b}$.
We recall further that the binary decision
between each pair of coordinates is performed via
a anti-symmetric function $f : \bits^n \to \bits$ so that
$f(-x) = -f(x)$ for all $x \in \bits^n$.
Finally, the tournament $G_K = G_K(\sigma; f)$ is defined
by having $(a > b) \in G_K$ if and only if $f(x^{a > b}) = 1$.\\

In order to obtain an upper bound on the probability that $1$ is the
unique maximum, define
$f^{a,b}(\sigma_1,\ldots,\sigma_n) = (1+f(x^{a > b}))/2$. Then the probability
that $1$ is the maximum is given by:
\[
\E[\prod_{a=2}^k f^{1,a}].
\]
Using~(\ref{eq:corr_bd}) we obtain that
\[
\E[\prod_{a=2}^k f^{1,a}] \leq
\StabThrmax{(\frac{1}{\sqrt{3}},\ldots,\frac{1}{\sqrt{3}})}{\frac{1}{2},\ldots,\frac{1}{2}} + \eps,
\]
where $\eps \to 0$ as $\tau \to 0$; there are $k-1$
expected values all given by $1/2$ as $\E[f^{1,a}] = 1/2$ for all $a$;
the $k-2$ values of $\rho$ all bounded by $1/\sqrt{3}$ by
Lemma~\ref{lem:rho_max_vote}.

By Lemma~\ref{lem:rho_bound} we now obtain
\[
\E[\prod_{a=2}^k f^{1,a}] \leq k^{-2 + o(1)}.
\]
Taking the union bound on the $k$ possible maximal values we obtain that
in Condorcet voting the probability that there is a unique max is at most
$k^{-1+o(1)}$ as needed.
\end{proof}

\subsection{Majority and Tightness} \label{sec:maj}
The majority function shows the tightness of Theorem~\ref{thm:MIMP} and
Theorem~\ref{thm:mibfc}.

\subsubsection{Tightness for Prediction}
For the tightness of Theorem~\ref{thm:MIMP}, write:
\[
X = \frac{\sum_{i=1}^n x_i}{\sqrt{n}}, \quad
Y = \frac{\sum_{i=1}^n x_i y_i}{\sqrt{\rho n}},
\]
where $x_i$ is the intended vote of voter $i$ and $y_i = 1$ if voter $i$ was queried, $y_i = 0$ otherwise. $X$ is the total bias of the actual vote and
$Y$ the bias of the sample.

Note that $(X,Y)$ is asymptotically a normal vector with covariance matrix
\[
\left( \begin{smallmatrix}
1 & \sqrt{\rho} \\
\sqrt{\rho} & 1
\end{smallmatrix} \right)
\]
Therefore asymptotically we obtain
\[
\E[\sgn(X) \sgn(Y)] = 2 \P[X >0, Y > 0] - 1 = 2 \StabThr{\rho}{1/2,1/2} - 1 =
{\textstyle \frac{2}{\pi}} \arcsin \sqrt{\rho}.
\]

\subsubsection{Tightness in Condorcet Voting}

\noindent
The tightness in Theorem~~\ref{thm:mibfc} follows from~\cite{RinottRotar:01}
and~\cite{NiemiWeisberg:68}. We briefly sketch the main steps of the proof.

For each $a$ and $b$ let
\[
X^{a > b} = \frac{\sum_{i=1}^n x_i^{a > b}}{\sqrt{n}},
\]
be the bias preference in a majority vote towards $a$.
By the CLT, all the random variables $X^{a > b}$ are asymptotically $N(0,1)$.

Consider the random variables $X^{1 > 2},\ldots,X^{1 > k}$.
Note that this set of variables is exchangeable (they are identically
distributed under any permutation of their order). Moreover,
\[
\E[X^{1 > 2} X^{1 > 3}] =
\frac{1}{n} \sum_{i=1}^n \E[x^{1 > 2}_i x^{1 > 3}_i] = 1/3.
\]
By the CLT the limiting value (as $n$ to $\infty$) of
\[
\lim_{n > \infty} \P[X^{1 > 2} > 0, \ldots, X^{1 > k} > 0] =
\P[N^{1 > 2} > 0, \ldots, N^{1 > k} > 0]
\]
where $(N^{1 > a})$
is an exchangeable collection of normal $N(0,1)$ random variables,
the correlation between each pair of which is $1/3$.
The results of~\cite{RinottRotar:01} imply that as $k \to \infty$:
\[
\P[N^{1 > 2} > 0, \ldots, N^{1 > k} > 0] \sim
2 \Gamma(2) \sqrt{2 \pi} \frac{\sqrt{\ln k}}{k^2},
\]
This in turn implies that the probability of a unique max for majority
voting for large $k$ as $n \to \infty$ is given by:
\[
(2 + o(1)) \Gamma(2) \sqrt{2 \pi} \frac{\sqrt{\ln k}}{k},
\]
showing the tightness of the result up to sub-polynomial terms.

\subsubsection{The probability that majority will result in linear order}
Here we prove Proposition~\ref{prop:maj_linear} and show
that the probability that majority will result in a linear order is
$\exp(-\Theta(k^{5/3}))$. We find this asymptotic behavior quite surprising.
Indeed, given the previous results that the probability that there is a unique max is $k^{-1+o(1)}$, one may expect that the probability that the order
is linear would be
\[
k^{-1+o(1)} (k-1)^{-1+o(1)} \ldots = (k!)^{-1+o(1)}.
\]
However, it turns out that there is a strong negative correlation between the
event that there is a unique maximum among the $k$ candidates and that among
the other candidates there is a unique max.

\begin{proof}
We use the multi-dimensional CLT. Let
\[
X_{a > b} = \frac{1}{\sqrt{n}} \left(
|\{\sigma : \sigma(a) > \sigma(b)\}| - |\{\sigma : \sigma(b) > \sigma(a) \}|
\right)
\]
By the CLT at the limit the collection of variables $(X_{a > b})_{a \neq b}$
converges to a joint Gaussian vector $(N_{a > b})_{a \neq b}$ satisfying
for all distinct $a,b,c,d$:
\[
N_{a > b} = - N_{b > a}, \quad
\Cov[N_{a > b}, N_{a > c}] = \frac{1}{3}, \quad
\Cov[N_{a > b}, N_{c > d}] = 0.
\]
and $N_{a > b} \sim N(0,1)$ for all $a$ and $b$.

We are interested in providing bounds on
\[
\P[\forall a > b: N_{a > b} > 0]
\]
as the probability that the resulting tournament is an order is obtained
by multiplying by a $k! = \exp(\Theta(k \log k))$ factor.

We claim that there exist independent $N(0,1)$ random variables $X_a$ for $1 \leq a \leq k$ and
$Z_{a > b}$ for $1 \leq a \neq b \leq k$ such that
\[
N_{a > b} = \frac{1}{\sqrt{3}} (X_a - X_b + Z_{a > b})
\]
(where $Z_{a > b} = - Z_{b > a}$).
This follows from the fact that the joint distribution of Gaussian random
variables is determined by the covariance matrix (this is noted
in the literature in~\cite{NiemiWeisberg:68}).

We now prove the upper bound. Let $\alpha$ be a constant to be chosen later.
Note that for all $\alpha$ and large enough $k$ it holds that:
\[
\P[|X_a| > k^{\alpha}] \leq \exp(-\Omega(k^{2 \alpha})).
\]
Therefore the probability that for at least half of the $a$'s in the interval
$[k/2,k]$ it holds that $|X_a| > k^{\alpha}$ is at most
\[
\exp(-\Theta(k^{1 + 2 \alpha})).
\]

Let's assume that at least half of the $a$'s in the interval $[k/2,k]$ satisfy
that $|X_a| < k^{\alpha}$. We claim that in this case the number
$H_{k/4}[-k^{\alpha},k^{\alpha}]$
of pairs $a > b$ such that
$X_a, X_b \in -[k^{\alpha},k^{\alpha}]$ and $X_a - X_b < 1$ is $\Omega(k^{2-\alpha})$.

For the last claim partition the interval
$[-k^{\alpha},k^{\alpha}]$ into sub-intervals of
length $1$ and note that at least $\Omega(k)$ of the points belong to
sub-intervals which contain at least $\Omega(k^{1-\alpha})$ points.
This implies that the number of pairs $a > b$ satisfying
$|X_a - X_b| < 1$ is $\Omega(k^{2-\alpha})$.

Note that for such pair $a > b$ in order that $N_{a > b} > 0$ we need
that $Z_{a>b} > -1$ which happens with constant probability.

We conclude that given that half of the $X$'s fall in $[-k^{\alpha},k^{\alpha}]$
the probability of a linear order is bounded by
\[
\exp(-\Omega(k^{2-\alpha})).
\]
Thus overall we have bounded the probability by
\[
\exp(-\Omega(k^{1 + 2 \alpha})) + \exp(-\Omega(k^{2-\alpha})).
\]
The optimal exponent is $\alpha=1/3$ giving the desired upper bound.

For the lower bound we condition on $X_a$ taking value in
$(a,a+1) k^{-2/3}$. Each probability is at least $\exp(-O(k^{2/3}))$ and therefore the probability that all $X_a$ take
such values is
\[
\exp(-O(k^{5/3})).
\]
Moreover, conditioned on $X_a$ taking such values the probability that
\[
Z_{a > b} > X_b - X_a,
\]
for all $a > b$ is at least
\[
\left( \prod_{i=0}^{k-1} \Phi(i)^{k^{2/3}} \right)^k \geq \left( \prod_{i=0}^{\infty} \Phi(i) \right)^{k^{5/3}} =
\exp(-O(k^{5/3})).
\]
This proves the required result.
\end{proof}

\section{Applications to Hyper-Graphs and Additive Properties} \label{sec:graphs}
In this section we prove Theorem~\ref{thm:ARIT} and give a few examples. The basic idea in the applications presented so far was that in order to
bound correlation between $k$ events of low influences, it suffices to know how to bound the correlation between the first $k-1$ and the last one.
For low influence events, using the invariance principle,
one obtains bounds coming from half spaces in Gaussian space,
or majority functions in the discrete space.

The applications presented now will be of different nature.
We will be interested again in
correlation between $k$ events -- however, we will restrict to correlation
measures defined in such a way that any pair of events are un-correlated.
While this is a much more restrictive setting, it allows one to
obtain exact results and not just bounds. In other words, we obtain that such
correlation measures for low-influence events depend only on the measure
of the sets but not on any additional structure. While this may sound
surprising, it in fact follows directly the invariance principle together with
the fact that for jointly Gaussian random variables,
pairwise independence implies independence. We first prove
Theorem~\ref{thm:ARIT}.

\begin{proof}
The proof follow immediately from Theorem~\ref{thm:MIST_easy},
Proposition~\ref{prop:MIST_easy} and Lemma~\ref{lem:cheeger}.
\end{proof}.

\begin{example}
Consider the group $\Z_m$ for $m > 2$. We will be interested in linear
relations over $\Z_m^k$. A {\em linear relation} over $\Z_m^k$ is given
by $L = (L_0,\ell_1,\ldots,\ell_k)$
such that $\ell_i \neq 0$ for all $i \geq 1$ and $\ell_i$ and $m$ are
co-prime for all $i \geq 1$.  We will restrict to the case $k \geq 3$.
We will write $L(x)$ to denote the logical
statement that $\sum_{i=1}^k x_i \ell_i \mod m \in L_0$.
Given a set $A \subset \Z_m^n$ we will denote
\[
L(A^k) = |\{ x \in A^k : L(x) \}|,
\]
and $\mu_L$ the uniform measure on $L(A^k)$.
We note that for every linear relation we have that $\mu_L$ is pairwise smooth
and that if the set $L_0$ is of size at least $2$ then $R$ is connected.

We now apply Theorem~\ref{thm:ARIT} to conclude that for low
influence sets $A \subset \Z_m^n$, the number of $k$ tuples
$(x_1,\ldots,x_k) \in A^k$ satisfying $\sum x_i \mod m \in L_0^n$ is
\begin{equation} \label{eq:ind_A}
(1 \pm o(1)) |A|^k \left( \frac{L_0}{m} \right)^n.
\end{equation}

For general sets $A$ we conclude that we have
$\underline{A}^S \subset A \subset \overline{A}^T$ where
$|S| = O(1), T = O(1)$ and~(\ref{eq:ind_A}) holds for both $\underline{A}^S$
and $\overline{A}^T$.
\end{example}

\begin{example}
We may consider much more general relations. For example we may take
$P$ be a polynomial in $x_1,\ldots,x_r \in \Z_r$ and $Q$ to be a polynomial
in $x_{r+1},\ldots,x_k$ such that $P$ and $Q$ both have roots.
Then we can look at the relation defined by the zeros of $PQ$.
It is easy to check that $R$ is connected. Therefore
 it follows that if $A \subset Z_m^n$
is set all of whose influences are lower than $\tau$ then
\[
c_1(\tau,\frac{|A|}{|\Z_m^n|}) |R^n|
\leq |A^k \cap R^n| \leq
c_2(\tau,\frac{|A|}{|\Z_m^n|}) |R^n|
\]
where $c_1$ and $c_2$ are two positive functions. Again if $A$ is not
of low influences then there exist finite sets of coordinates $S$ and $T$
such that $\underline{A}^S \subset A \subset \overline{A}^T$ and the conclusion
holds for $S$ and $T$.
\end{example}

\bibliographystyle{abbrv}
\bibliography{all,my}

\end{document}